\documentclass{article}

\title{$k$-linear Morita theory}
\author{Matteo Doni\footnote{Università degli Studi di Milano, Milan. Email address: matteo.doni@unimi.it}}
\input{CategoricalArticle.sty}
\begin{document}
\newlength{\myindent} 
\setlength{\myindent}{\parindent}
\parindent 0em 
\maketitle
\begin{abstract}
In this paper, we prove the standard comparison used by mathematicians between the idempotent complete pretriangulated dg-categories, over a unitary and commutative ring $k$, and the idempotent complete $k$-linear stable $\infty$-categories. Our approach is completely included in the $\infty$-categorical theory. To achieve the target we will reinterpret the Morita theory for dg-categories and we set the Morita theory for $k$-linear stable $\infty$-category. 

\end{abstract}
\setcounter{tocdepth}{1}
\tableofcontents

\section*{Introduction}
Since their first appearance in the 60's by Verdier the theory of triangulated category has gained interest. Today triangulated categories are central objects in
many fields of mathematics, ranging from homotopy theory to algebraic and symplectic
geometry. The Russian school played a central role in this development, particularly due to the efforts of Bondal, Orlov, and Kapranov. According to
Kontsevich, the most fundamental statement of mirror symmetry is as an equivalence of
categories between the derived category of coherent sheaves on a Calabi-Yau variety and
the Fukaya category of its mirror. The derived category $D(\mathcal{A})$ of an abelian category $\mathcal{A}$ is certainly a very intriguing and ubiquitous mathematical object. It is well known that $D(\mathcal{A})$ is a triangulated category and this leads to some complications: for example, the famous \enquote{non-functoriality of the cone} problem.
The disadvantage is that during the transition to the derived category $D(\mathcal{A})$ from the abelian $\mathcal{A}$ some \enquote{higher} information is lost.
The strategy to solve this issue is to lift $\mathcal{A}$ to a higher category $\mathcal{C}$, whose homotopy category $h\mathcal{C}$ is triangulated and such that there exists an equivalence of triangulated categories $h\mathcal{C}\simeq D(\mathcal{A})$: dg-category (over a suitable ring) and $\infty$-category theories are two of the common higher category theories used to define $\mathcal{C}$. In this paper, we are interested in them. 

Not all homotopy categories of a dg-category or of an $\infty$-category are triangulated, but when they are, they are called pretriangulated and stable, respectively.
We must apply a first correction in our previous statement: in this paper, we are interested in pretriangulated dg-categories and stable $\infty$-category theories. 
    
Algebraic geometers are focused on studying the existence and uniqueness of such lifts.
An in-depth study of this issue in a pretriangulated dg-categorical setting was made by V. Lunts, D. Orlov A. Canonaco, A. Neeman, P. Stellari, see \cite{canonaco2017tour}\cite{canonaco2018localizations}\cite{canonaco2022uniqueness}\cite{lunts2010uniqueness}; instead B. Antieua tackles it in stable $\infty$-categorical style, see \cite{antieau2018uniqueness}.
Actually, the \enquote{lift problem} is a bit general, and it can be formulated as follows: \textit{Let $T$ be a triangulated category, does a dg-categorical (or a stable $\infty$-categorical) lift exist? And if so, is it unique?}

A triangulated category $T$ is called \textit{algebraic} if there exists a pretriangulated dg-categorical lift and, it is called \textit{topological} if a stable $\infty$-categorical lift exists.
Let $\mathcal{A}$ be a Grothendieck abelian category. The derived $\infty$-category $\mathcal{D}(\mathcal{A})$ exists and it is a stable $\infty$-categorical lift of the ordinary derived category $D(\mathcal{A})$. The full dg-subcategory of h-injective objects $\mathcal{C}_{\text{dg}}(h$-$\text{inj})\subseteq\mathcal{C}_{\text{dg}}(\mathcal{A})$ is a dg-categorical lift of $D(A)$, \cite[Example 3.4]{canonaco2017tour} and \cite[Proposition 1.3.5.15]{HA}. 
Recall that $I\in K(\A)$ is called h-injective (or K-injective, as in [62]) if $\mathrm{Hom}_{K(\A)}(A, I)$ = 0 for every acyclic complex $A\in K(\A)$.

This example is useful because it allows us to understand the difference between topological and algebraic triangulated categories. First of all, given an abelian category $\mathcal{A}$, it is interesting that the idea of finding $\mathcal{D}(A)$ and $\mathcal{C}_{\text{dg}}(h$-$\text{inj})$ is the same one: it is just a translation from one theory to the other using the Lurie differential graded nerve $N_{dg}(-)$, \cite[Construction 1.3.1.6]{HA}. 
This strategy is general and not specific to this example. Any dg-categorical lift can be translated into a stable $\infty$-categorical lift via $N_{dg}(-)$, as shown in \cite[Proposition 1.3.1.10 and Remark 1.3.1.11]{HA}. Consequently, every algebraic triangulated category is also topological. Shortly, we will demonstrate that the converse does not hold.

The second conclusion follows from a specific case. Let $X$ be a scheme over $k$, it is well known that the category $QCoh(X)$ of quasi-coherent sheaves over $X$ is abelian, then the derived category $D(QCoh(X))$ is a topological algebraic triangulated category.
The dg-categorical lift $\mathcal{C}_{\text{dg}}(h$-$\text{inj})$ is naturally a $k$-linear category, since it is dg-category over $k$; instead, a priori, the stable $\infty$-categorical lift has no natural $k$-linear structure.
Every scheme $X$, in Grothendieck's conception, is at least over $\mathbb{Z}$, then each dg-categorical lift is at least $\mathbb{Z}$-linear: that is why for algebraic geometers the interesting lifts are those with a linear structure.
Of course, one can ask if there is a canonical way to give a $k$-linear structure, for a suitable ring $k$, to a stable $\infty$-categorical lift. Generally, this is not possible.
The category $Sp^{\Sigma}$ of symmetric spectrum is topological, its stable $\infty$-categorical lift is $\Sp$, but it is not algebraic, \cite[Proposition 3.6]{canonaco2017tour}.
So topological triangulated categories are more than algebraic ones: the problem is that there is not a canonical linearization of $Sp^{\Sigma}$. 

One can also ask the opposite question: are there stable $\infty$-categories with a (canonical) linearization?

Our first example comes to our aid. Since $\mathcal{D}(A)$ is the image of $\mathcal{C}_{\text{dg}}(h$-$\text{inj})$ via the $\infty$-functor $N_{dg}(-)$ it has a linearization. It shows that some stable $\infty$-categories have a canonical linearization. When this happens they are called $k$-linear stable $\infty$-categories.


An important achievement in this regard is \cite[Meta Theorem]{antieau2018uniqueness}. B. Antieau proves that the presence of a separated or anticomplete $0$-cosimplicial t-structure guarantees the existence of a canonical $\Z$-linear enrichment. So, for triangulated categories admitting a certain structure, stable $\infty$-categorical lifts have canonical $\Z$-linearization.

One point needs to be made before you can define $\mathcal{D}(\mathcal{A})$, for what we know, you first have to define a dg-category. The use of the dg-category theory to define the stable $\infty$-category $\D(\mathcal{A})$ is not accidental, but it is because this theory is more traditional and studied when working in a linear setting. 

The category of symmetric spectra $Sp^{\Sigma}$ is a fundamental category for algebraic topologists, so they are concerned with a wider range of lifts. So, the term $k$-linear means that the mathematical object comes from geometry. To highlight this fact, we add the term \enquote{not-geometrical} when dealing with a concept without considering $k$-linearity

The algebraic topologist point of view leads to a spectral algebraic geometry theory, see \cite{SS}; the algebraic geometric point of view leads to derived algebraic geometry theory, see \cite{toen2014derived}.

Our Main result, \Cref{thMaintheorem}, is a meeting point between them. It states that:

\begin{theorem}[Main Theorem]
\label{MainTheoremIntroduzione}
Let $k$ be a discrete $\E$-ring. Let $\Cat^{\D(k),perf}$ be the $\infty$-category of idempotent complete pretriangulated dg-categories over $k$, and let $LinCat^{perf}_{\Hk}$ the $\infty$-category of idempotent complete $k$-linear stable $\infty$-categories. There exists an equivalence of $\infty$-categories:
\begin{equation}
\Cat^{\D(k),perf}\simeq LinCat^{perf}_{\Hk}.
\end{equation}
\end{theorem}

It is always useful to work with several theories on the same problem; we call both $k$-linear Morita theory. The two settings in our Main Theorem have different features. Dg-categories are the simplest and most hands-on enhancement
of triangulated categories.
Lurie develops the second formalism in a series of enormous and instrument-filled books. The theory of $k$-linear stable $\infty$-categories is more flexible, as the theory of localization proves,
and the standard enhancement of triangulated categories in many mathematical areas, such as homotopy theory.
A dictionary between different enhancements makes it
possible, in any given situation, to adopt the formalism which is most convenient; while
at the same time allowing for the transfer of results and techniques that are native to
the other formalism. As I mentioned, these two different enhancements have actually
quite different features, so being able to switch freely between them is particularly valuable.

In \Cref{MainTheoremIntroduzione}, the condition of idempotent completeness appears, which we have not yet discussed. We remark that this property serves to handling size issues. Indeed, Lurie shows that this notion, in the context of $\infty$-categories, has at least two distinct incarnations: one categorical, as in \Cref{exIdempo}, and another related to the compact $\infty$-subcategory of a presentable category, \eqref{eqZZoCA}, which is an essentially small $\infty$-category.

The fluidity of this notion is crucial in this paper because it allows us to define a symmetric monoidal structure for the $\infty$-categories where the non-geometrical Morita theory is embodied, $\Cat^{perf}$ and $Pr_{\omega}^{L, St}$, see \Cref{subsecCatperfPRomega}.
$\\$


\section*{Related works}
This paper is the last of a series \cite{DoniCategorical,DoniHigherCategorical}, aimed at developing and comparing the Morita theory for dg-categories, over a commutative unitary ring $k$, and $k$-linear $\infty$-categories in a $\infty$-categorical setting.
Here we finally set this theory and we call it the \textit{$k$-linear Morita theory.}
The $\infty$-categorical Morita theory, without considering that $k$-linearity, was set up in the pioneering article \cite{blumberg2013universal} by A.J. Blumberg, D. Gepner, and G. Tabuada. In this paper, we call it the non-geometrical Morita theory.
$\\$

Our Main Theorem, \Cref{thMaintheorem}, was a folklore theorem until 2013; in fact, in 2010 D. Ben-Zvi, J. Francis, and D. Nadler in \cite{IntegralTransf} and in 2012 D. Gaitsgory in \cite{Gaitsgory} used it, although it had not yet been proven.

In 2013 L. Cohn proved it in a preprint, and in 2016, he updated the paper to its latest version. The preprint has been improved along with the young and exuberant $\infty$-category and spectral algebraic geometry theories.   

It is also interesting to note that J. Lurie conjectured our Main Theorem in our version, \cite[Remark 1.1.0.2]{HA}, and in the additive version, \cite[D.1.2.3]{SAG}.
    
The $\infty$-category theory and the spectral algebraic geometry theory have made further progress. Therefore, we think that a reproof can be useful to make this result clearer for more mathematicians.
An important difference between our proof and Cohn's is that our approach is completely $\infty$-categorical.



\section*{Structure of the work}
This paper is composed of six sections.
In \Cref{secAddingColimits}, we recall how in $\infty$-category theory it is possible to add (chosen) colimits.   
In \Cref{subsecStablepresentable}, we review key results, definitions, and construction from higher algebra. 
For a reader used to the work of Lurie, \cite{HTT, HA}, the first two sections can be skipped.

In \Cref{subsecCatperfPRomega}, we review the non-geometrical Morita theory constructed in \cite{blumberg2013universal}. Besides establishing definitions for subsequent sections, the primary objective of this chapter is to demonstrate that the categories $\Cat^{perf}$ and $Pr^{L,St}_{\omega}$ possess certain higher algebraic structures (\Cref{rmkSptensorizationonCatPerfPrSt}). Some of these structures were previously identified, but their proof was either flawed or unconvincing.
With this section the prerequisite part concludes, paving the way for the substantive sections of the paper where we tackle for the first time the problem of defining the theory of $k$-linearization using enriched $\infty$-categories (Lurie has previously addressed this topic but with a different setting \cite[\S D]{SAG}).

In \Cref{secDg}, we reinterpret the Morita theory for dg-categories, \cite{toen2007homotopy, tabuada2005structure}, using the theory of enriched $\infty$-category. This theory is classical in derived algebraic geometry but our $\infty$-categorical approach is not: this is a substantial difference between our proof and Cohn's.

In \Cref{secklinearstableCat}, we investigate linearization from the viewpoint of spectral algebraic geometry (or algebraic topology). We provide a conceptual overview of linearization and define Morita theory for $k$-linear $\infty$-categories. In this part, we will come to realize that there exists an incredible number of equivalent definitions of idempotent complete $k$-linear stable $\infty$-categories, and the main outcome of this section is to demonstrate their equality (\Cref{propNewDefinitionProj}).

Finally, in \Cref{secRacks}, we present our principal result, \Cref{thMaintheorem}. It is valuable that we define when a dg-functor is a Morita equivalence, \Cref{defMoritaEquivalence}, and find the ($k$-linear) Morita localization, \Cref{propklinearLocalMorita}.

\section*{Acknowledgement}
I extend my gratitude to Paolo Stellari for his guidance throughout the recent years, which has enabled me to consolidate my ideas.

\section*{Notation}
\begin{itemize}
   


    \item We adhere to Lurie's homological convention, as described in \cite{HA};

    


    \item  with a slight abuse of notation, we denote by $\pi_{0}:\mathcal{S}_{*}\to Set$ the composition of $\infty$-functors $U\pi_{0}$ where $\pi_{0}:\mathcal{S}\to Set$ is the component-of-space $\infty$-functor,\cite[Example 3.6.2]{RiehlCTIC}, and $U:\mathcal{S}_{*}\to\mathcal{S}$ is the canonical forgetful $\infty$-functor;


    \item unless otherwise indicated by the symbol $\dashv$, the following conventions will be used: in a horizontally written adjunction, the left (right) adjoint has the arrowhead pointing to the right (left);




    \item we use the symbol $\mapsto$ to indicate where a functor sends objects or where an $\infty$-functor, up to equivalence, maps objects. For example, $-+1:\Z\to\Z:n\mapsto n+1$.

    \item let $\C$ be an $\infty$-category. By $h\C$, we denote the ordinary category whose objects are the same as those of $\C$ and whose hom-set between two objects $x, y$ is the set $\pi_0(\C(x, y))$. We call $h\C$ the \textit{homotopy category of $\C$}.
    For a comprehensive presentation of $h\C$, we recommend reading \cite[\S 1.1]{RiehlElements}.

    \item throughout this paper, we fix $k$ a commutative ring with unit (or a discrete $\E$-ring);

    \item it is very likely that if a notation is used without defining it, you can find its definition in the classical references \cite{HA, HTT} or that we set it in \cite{DoniHigherCategorical}.





\end{itemize}

\section{Recall of higher categories: adding colimits}
\label{secAddingColimits}

In this part, we recall how to add colimits to an $\infty$-category and how this technique is functorial and monoidal. For a detailed presentation, see \cite[\S 4.8]{HA}.

\begin{notation}
We denote by $\widehat{\Cat}$ the $\infty$-category of non-necessarily small $\infty$-categories.
\end{notation}

\begin{definition}
\label{defCatK}
Let $\mathscr{K}$ be a collection of simplicial sets. We let $\Cat(\mathscr{K})$ ($\widehat{\Cat}(\mathscr{K})$) denote the $\infty$-subcategory of $\Cat$ ($\widehat{\Cat}$) spanned by those small (non-necessarily small) $\infty$-categories $\C$ which admit $K$-indexed colimits for each $K \in \mathscr{K}$, and whose morphisms are $\infty$-functor $f : C \to D$ which preserve $K$-indexed colimits for each $K \in\mathscr{K}$.
\end{definition}

We need the following proposition, which is at the same time a definition.

\begin{proposition}[{\cite[Proposition 5.3.6.2]{HTT}}]
\label{defPKH}

Let $\mathscr{K}$ be a collection of simplicial sets, $\C$ an $\infty$-category, and $\mathscr{R}=\{d_{\alpha}:\mathbb{1}\to \C^{K_{\alpha}^{\rhd}}\}$ a collection of diagrams in $\C$. Assume that each $K_{\alpha}$ belongs to $\mathscr{K}$. Then there exists a new $\infty$-category $P_{\mathscr{R}}^{\mathscr{K}}(\C)$ and a $j:\C\to P_{\mathscr{R}}^{\mathscr{K}}(\C) $ with the following properties:
\begin{itemize}
    \item[(1)] the $\infty$-category $P_{\mathscr{R}}^{\mathscr{K}}(\C) $ admits $\mathscr{K}$-indexed colimits;
    \item[(2)] For every $\infty$-category $\D$ which admits $\mathscr{K}$-indexed colimits, composition with $j$ induces an equivalence of $\infty$-categories \[(-)\circ j:Fun_{\mathscr{K}}(P_{\mathscr{R}}^{\mathscr{K}}(\C),D)\to Fun_{\mathscr{R}}(\C,\D):Lan_{j}(-);\]
\end{itemize}
Moreover, if every member of $\mathscr{R}$ is already a colimits diagram in $\C$, then we have in addition:
\begin{itemize}
    \item[(3)] the $\infty$-functor $j$ is fully faithful. 
\end{itemize}
\begin{notation}
    If $\mathscr{K}$ or $\mathscr{R}$ is empty we will erase the apex or respectively the subscript from $P_{\mathscr{R}}^{\mathscr{K}}(\C)$.
\end{notation}
\end{proposition}


Let $\K$ be a collection of simplicial sets. In \cite[Corollary 4.8.1.4]{HA}, the author proves that $\Cat(\K)$ ($\widehat{\Cat}(\K)$) admits a symmetric monoidal structure and that the inclusion $\Cat(\K)\to \Cat$ ($\widehat{\Cat}(\K)\to \widehat{\Cat}$) is a lax symmetric monoidal $\infty$-functor. We denote the tensor of this category by $\otimes^{\K}$, and the unit is the $\infty$-category $P^{\K}(*)$.

Furthermore, if $\K \subseteq \mathscr{H}$ are collections of simplicial sets, then we obtain a monoidal $\infty$-functor from:
\begin{equation}
\label{localizationCatKH}
P^{\mathscr{H}}_{\K}(-):\Cat(\K) \to \Cat(\mathscr{H}).
\end{equation}

The following examples are important for us.

\begin{example}
\label{exPkhSmallColimits}
Let $\K$ be the collection of all small simplicial sets, then $P(\C)\simeq P^{\K}(\C)$. 
\end{example}

\begin{example}[{\cite[Corollary 4.4.5.15]{HTT}}]
\label{exIdempo} 
Let $Idem =N(R)$ be the nerve of the category $R$ with only one object $X$ and whose arrows are two, $e$ and the identity $id$, and such that, without considering the obvious cases, the composition of arrows is described by the equation $e \circ e = e$:
\begin{equation}
\label{eqIdem}
\begin{tikzcd}
	X & X & X.
	\arrow["e", from=1-1, to=1-2]
	\arrow["e", from=1-2, to=1-3]
	\arrow[""{name=0, anchor=center, inner sep=0}, "e"', bend right, from=1-1, to=1-3]
	\arrow["\circlearrowright"{description}, draw=none, from=1-2, to=0]
\end{tikzcd}
\end{equation}

Let $N=\{Idem\}$ be the collection of simplicial sets composed only of $Idem$, then \[P^{\mathcal{K}}(\C)\simeq Idem(\C):=Fun(Idem,X)\] is the idempotent completion of the $\infty$-category $\C$.

It is important to note that the morphisms in $\Cat(N)$ are $\infty$-functors between $\infty$-categories with no other properties: for a diagram with source an $Idem$, having a colimit is a diagram condition, see \cite[Corollary 4.4.5.14]{HTT}, so every $\infty$-functor preserves its colimit; in short, they are absolute colimits. We denote by $\Cat^{Idem}$ the $\infty$-category $\Cat(N)$ and we call its objects \textit{idempotent complete $\infty$-categories}.

The simplicial set $Idem$ is not a finite simplicial set but it is a $\alpha$-filtered simplicial set for each regular cardinal $\alpha$, \cite[Example 5.3.1.9]{HTT}.
\end{example}

\begin{example}
\label{exInd}
Let $\alpha$ be a regular cardinal and let $K_{\alpha}$ be the collection of all small $\alpha$-filtered simplicial sets, then $P^{K_{\alpha}}(\C)\simeq Ind_{\alpha}(\C)$. So, there is an $\infty$-functor \[Ind_{\alpha}(-):\Cat\to \Cat(K_{\alpha}).\]
The objects of a $\infty$-category $\C$ in $\Cat(K_{\alpha})$ are usually called \textit{Ind-objects}.

As usual, when $\alpha=\omega$ we write $Ind(-)$ instead of $Ind_{\omega}(-)$.  
\end{example}

\begin{example}
\label{exfinitedimpli}
Let $\F$ be the collection of finite simplicial sets. For a definition of finite simplicial set see \cite[Definition 3.6.1.1]{kerodon}.
$\Cat(\F)$ is the $\infty$-category of $\infty$-categories admitting finite colimits and whose morphisms are finite-colimits-preserving $\infty$-functors; as usually we denote $\Cat(\F)$ by $\Cat^{rex}$.

We know that $\Cat^{rex}$ has a symmetric monoidal structure whose unit is the $\infty$-category $P^{K}(*)\simeq \Top^{fin}$. We denote its tensor by $\otimes^{rex}$.
\end{example}

\begin{example}
\label{exCatRIdem}
Let $\F$ be the collection of finite simplicial sets and $N$ as in \Cref{exIdempo} we denote by $\Cat^{RIdem}$ the $\infty$-category $\Cat^{K\cup \F}$ and we call its objects \textit{idempotent complete cocomplete $\infty$-categories}.   
\end{example}

\begin{remark}
\label{rmkCatIdemsymmetricMonoidal}
Let $N$ and $\F$ be as in \Cref{exfinitedimpli} and \Cref{exCatRIdem} respectively. Both the $\infty$-categories $\Cat^{Idem}$ and $\Cat^{RIdem}$ have a symmetric monoidal structure. Furthermore, there exists a monoidal adjunction
\[\begin{tikzcd}[ampersand replacement=\&,sep=small]
	{\Cat^{Idem}} \&\&\& {\Cat^{RIdem}.}
	\arrow[""{name=0, anchor=center, inner sep=0}, "{P^{F}_{N}(-)}", shift left=3, from=1-1, to=1-4]
	\arrow[""{name=1, anchor=center, inner sep=0}, shift left=3, from=1-4, to=1-1]
	\arrow["\dashv"{anchor=center, rotate=-90}, draw=none, from=0, to=1]
\end{tikzcd}\]
We do not think that the above adjunction is a localization. 
\end{remark}

\begin{remark}
\label{rmkLocalizationrexRIdem}
Let $N$ be as in \Cref{exIdempo} and $\F$ be as in \Cref{exfinitedimpli}. Since every $\infty$-functor preserves idempotent object $\Cat^{RIdem}\subseteq \Cat^{rex}$ is a full $\infty$-subcategory, so the $\infty$-functor 
$P_{\F}^{N}(-):\Cat^{rex}\to \Cat^{RIdem}$ defines a localization: 
\[\begin{tikzcd}[ampersand replacement=\&,sep=small]
	{\Cat^{rex}} \&\&\& {\Cat^{RIdem}.}
	\arrow[""{name=0, anchor=center, inner sep=0}, "{P_{\F}^{N}(-)}", shift left=3, from=1-1, to=1-4]
	\arrow[""{name=1, anchor=center, inner sep=0}, shift left=3, hook', from=1-4, to=1-1]
	\arrow["\dashv"{anchor=center, rotate=-90}, draw=none, from=0, to=1]
\end{tikzcd}\]

\end{remark}

\begin{example}[{\cite[Example 5.3.6.8]{HTT}}]
\label{exTHeTrick}
Let $\alpha$ be a regular cardinal. Let $\K_{\alpha}$ denote the collection of all $\alpha$-small simplicial sets and let $\K$ be the class of all small simplicial sets. Let $\C$ be a small $\infty$-category which admits $\alpha$-small colimits. Then we have a canonical equivalence $P^{\K}_{\K_{\alpha}}(\C)\cong Ind_{\alpha}(\C)$. A simplicial set is $\omega$-small if it is finite, \cite[Example 5.4.4.4]{kerodon}, so $\K_{\omega}$ is equivalent to $\F$ in \Cref{exfinitedimpli}.
\end{example}

\begin{definition}[{\cite[Definition 5.4.2.1]{HTT}}] 
\label{defAccessible}
Let $\gamma$ be a regular cardinal. An $\infty$-category $\C$ is $\gamma$-accessible if there exists a small $\infty$-category $\C_0$ and an equivalence $Ind_{\gamma}(\C_0) \simeq \C$.
We will say that $\C$ is accessible if it is $\gamma$-accessible for some regular cardinal $\gamma$.
\end{definition}

\begin{definition}
\label{defcompactobjects}
Let $\C$ be an $\infty$-category which admits filtered colimits (for us it will always be an accessible $\infty$-category). An object $c\in \C$ is said to be compact if the corepresentable functor
\[\Hom_{\C}(c,-):\C\to \Top;\]
commutes with filtered colimits.
\end{definition}

\begin{notation}[{\cite[Remark 5.4.2.11]{HTT}}]
\label{RmkLEOLEOLEOLEO}
Let $\C$ be an $\infty$-category which admits filtered colimits (for us it will always be an accessible $\infty$-category) and $\gamma$ a regular cardinal.
We denote by $\C^{\gamma}$ the full $\infty$-subcategory spanned by the $\gamma$-compact objects of $\C$.
The $\infty$-category $\C^{\gamma}$ is essentially small.
\end{notation}

\section{Stable $\infty$-category and presentable $\infty$-categories}
\label{subsecStablepresentable}
\subsection{Stable and presentable $\infty$-category}
In \cite{DoniHigherCategorical}, we have already presented part of the theory of presentable and stable $\infty$-category. So if you find any notation or definition unclear, it is very likely that there is a more detailed explanation in the work just cited. 
Since this part is taken from \cite{HA, HTT}, these two books also contain a detailed presentation of these topics.
We start by recalling the definitions.

 There are a lot of definitions of stable $\infty$-category but we choose the most similar to the $dg$-categorical notion of pretriangulated $dg$-category.  

\begin{definition}[{\cite[Definitions 1.1.1.9, Corollary 1.4.2.27]{HA}}]
\label{defStable}
An $\infty$-category $\C$ is stable if it satisfies the following conditions:
\begin{itemize}
    \item There exists a zero object $\mathbf{0}$;
    \item Every morphism in $\C$ admits cofiber;
    \item The suspension $\infty$-functor \[\Sigma:\C\to \C:X\mapsto \operatorname{colim}(\mathbf{0}\leftarrow X\to \mathbf{0})\] is an equivalence of $\infty$-categories.
\end{itemize}
\end{definition}

\begin{example}
Recall in \Cref{exInd} we defined, for each regular cardinal $\gamma$, an $\infty$-functor \[Ind_{\gamma}(-):\Cat\to \Cat(K_{\alpha}).\]
\end{example}

\begin{definition}
\label{defPresentable}
an $\infty$-category $\C$ is presentable if it admits small colimits and it is accessible.
We say that an $\infty$-category $\C$ is \textit{compactly generated} if it is presentable and $\omega$-accessible.
\end{definition}

\begin{notation}
We use the usual notation $P_r^L\subseteq \widehat{Cat}_\infty$ to denote the $\infty$-subcategory spanned by presentable $\infty$-categories and whose morphisms are small-colimits-preserving functors.
We denote by $Pr^{L}_{cg}\subseteq Pr^{L}$ the full $\infty$-subcategory spanned by compactly generated $\infty$-categories.
\end{notation}
It is well known that there is a closed symmetric monoidal structure $(Pr^{L},\otimes^{L},\Top)$ such that the inclusion $\iota: Pr^{L}\to \widehat{\Cat}$ is lax monoidal and $\otimes^{L}$ preserves small colimits componentwise.  
We call the algebra objects of $Pr^{L}$ \textit{presentably monoidal $\infty$-categories}.






We need to know what a $t$-structure on a stable $\infty$-category is. Now we will recall its definition and explain some useful examples for us; for a complete introduction, see \cite[\textsection 1.2.1]{HA}.
Roughly speaking, giving a $t$-structure to a stable $\infty$-category $\C$ means adding additional information to $\C$ that allows $\C$ to be broken up into homologous parts, and by studying these parts it is possible to derive information on the whole of $\C$.

\begin{definition}[{\cite[Definition 1.2.1.1 Definition 1.2.1.4 ]{HA}}]
    Let $\C$ be a stable $\infty$-category. A triple $(\C,\C_{\geq 0},\C_{\leq 0})$ where $\C_{\geq 0}$ and $\C_{\leq 0}$ are two subcategories of $\C$ is a \textit{$t$-structure} for $\C$ if the triple $(h\C,h\C_{\geq 0},h\C_{\leq 0})$ is a $t$-structure for the triangulated category $h\C$.
  
\end{definition}

\begin{definition}[{\cite[Definition 1.2.1.11]{HA}}]
    Let $\C$ be a stable $\infty$-category equipped with a $t$-structure. \textit{The heart} $\C^{\heartsuit}$ of $\C$ is the full $\infty$-subcategory $\C_{\leq 0}\cap\C_{\geq 0}\subseteq\C $. For each $n\in\Z$, we let $\pi_{0}:\C\to \C^{\heartsuit}$ denote the $\infty$-functor $\tau_{\geq 0}\circ \tau_{\leq 0} \simeq \tau_{\leq 0} \circ\tau_{\geq 0}$, and let $\pi_{n}:\C\to\C^{\heartsuit}$ denote the composition of $\pi_{0}$ with the shift $\infty$-functor $X\mapsto X[-n]$. $\C^{\heartsuit}$ is an abelian (ordinary) category. 
\end{definition}

In the construction of the $\infty$-functor \eqref{eqChangeEnr}, it is important that the $\infty$-categories $\LModk$ and $\D(k)$ have a right complete $t$-structure.

\begin{definition}
\label{deftstructureComplete} 
    Let $\C$ be a stable $\infty$-category. A $t$-structure $(\C,\C_{\geq 0},\C_{\leq 0})$ is \textit{right complete} if $\C$ is equivalent to the direct limits of its right bounded part 
    \[
    \C\simeq \lim(\dots \xrightarrow{\tau_{\leq n-1}}\C_{\geq n-1}\xrightarrow{\tau_{\geq n}}\C_{\geq n}\to \dots).
    \]
\end{definition}

\begin{remark}
\label{rmkCompleteTstructeAndDualized}
We make two brief remarks about the last definition.
Note that in \Cref{deftstructureComplete}, we defined the right completeness as a limit of a diagram $D:\Z\to \Cat$ but in \Cref{eqD(k)LModK} as a limit of a diagram $H:\N^{op}\to\Cat$. Now we explain why the two limits have the same apex.  
There is a canonical adjunction: 
    \[ i:\N^{op}\rightleftarrows\Z: R;   \]
    such that, by straightforward computation, the equality $H=Di$ holds.
    Every left adjoint is an initial $\infty$-functor \cite[Proposition 2.4.6]{RiehlElements}, and so we conclude that $\lim(H)\cong \lim(Di)\cong \lim(D)$;

\end{remark}

Moreover, we need a method to compare $t$-structures.
\begin{definition}[{\cite[Definition 1.3.3.1]{HA}}] Let $\C$ and $\D$ be stable $\infty$-categories equipped with $t$-structures. We will say
that a $\infty$-functor $f : \C \to \D$ is \textit{right t-exact} if it is exact and carries $\C_{\geq 0}$ into $\D_{\geq 0}$. We say that $f$ is \textit{left t-exact} if it is exact and carries $C_{\leq 0}$ into $\D_{\leq 0}$.
\end{definition}

In the next part, we list the $t$-structures that are important to us.

\begin{example}
 \label{extstructureD(k)-}
 
Let $D(k)^-_{\leq 0}$ be the full $\infty$-subcategory of $\D(k)^-$ spanned by chain complexes with trivial positive homology groups, i.e., $H_{n}(C_*)\simeq 0$ if $ n>0 $.
Let $D(k)^-_{\geq 0}$ be the full $\infty$-subcategory of $\D(k)^-$ spanned by chain complexes with trivial negative homology groups, i.e., $H_{n}(C_*)\simeq 0$ if $n<0$.
 The triple $(D(k)^-,D(k)^-_{\leq 0},D(k)^-_{\geq 0})$ forms a $t$-structure.  
 In particular $\D(k)^{\heartsuit}\simeq k\text{-}Mod$.
\end{example}

\begin{example}[{\cite[Proposition 1.3.5.21]{HA}}]
\label{extstructureD(k)}
The $\infty$-category $\D(k)$ is compactly generated and stable.
Let $\D(k)_{\leq 0}$  be the full $\infty$-subcategory of $\D(k)$ spanned by chain complexes with trivial positive homology groups, i.e., $H_{n}(C_*)\simeq 0$ if $n>0 $. 
Let and $\D(k)_{\geq 0}$ be the full $\infty$-subcategory of $\D(k)$ spanned by chain complexes with trivial negative homology groups, i.e., $H_{n}(C_*)\simeq 0$ if $n<0 $.
 The triple $(D(k),D(k)_{\leq 0},D(k)_{\geq 0}) $ forms a right complete $t$-structure.  
 In particular $\D(k)^{\heartsuit}\simeq k\text{-}Mod$.
\end{example}

For us, there is another important stable $\infty$-category with a $t$-structure, but before we can consider it, we need a definition.


\begin{example}[{\cite[Proposition 7.1.1.5]{HA}\cite[Proposition 7.1.1.13]{HA}}]
\label{exSperiamoUltimoDaCitare}
The $\infty$-category $\LModk$ is compactly generated and stable.  
Let $ \LModk^{\leq 0}$ be the full $\infty$-subcategory of $\LModk$ spanned by modules with trivial positive homotopy groups, i.e., $\pi_{n}(M)\simeq 0$ if $n>0 $.
Let $\LModk^{\geq 0}$ be the full $\infty$-subcategory of $\LModk$ spanned by modules with trivial negative homotopy groups, i.e., $\pi_{n}(M)\simeq 0$ if $n<0 $.
The triple $(\LModk,\LModk^{\leq 0},\LModk^{\geq 0}) $ forms a right complete $t$-structure. 
In particular, $\LModk^{\heartsuit}\simeq k\text{-}Mod$.   
\end{example}


\begin{notation}
Let $\Cat^{Lex}\subseteq\Cat$ be the $\infty$-subcategory spanned by small $\infty$-categories with finite limits and whose morphisms preserve finite limits as $\infty$-functors.
Let $\Cat^{Ex}\subseteq\Cat$ be the $\infty$-subcategory spanned by small stable $\infty$-categories and whose morphisms are exact $\infty$-functors: an $\infty$-functor between stable $\infty$-categories is exact if it preserves finite colimits and finite limits. 
\end{notation}

The construction of the $\infty$-category of $\C$-spectra $Sp(\C)$ is left adjoint to two adjunctions \cite{HA,heine2023equivalence}:

\begin{equation}
\label{eqLocLexEx}
\begin{tikzcd}
	{\Cat^{Lex}} && {\Cat^{Ex};}
	\arrow[""{name=0, anchor=center, inner sep=0}, "{Sp(-)}", shift left=3, from=1-1, to=1-3]
	\arrow[""{name=1, anchor=center, inner sep=0}, "i", shift left=2, hook', from=1-3, to=1-1]
	\arrow["\dashv"{anchor=center, rotate=-90}, draw=none, from=0, to=1]
\end{tikzcd}
\end{equation}

For us, there is another important stable $\infty$-category with a $t$-structure, but before we can consider it, we need a definition.






Now, we study the $\infty$-categories which are both presentable and stable. Having these two properties is very interesting because it allows these categories to be seen as $\Sp$-modules objects. In other words, it is possible to translate a higher category concept into a higher algebraic concept.

\begin{notation}
 With $Pr^{L,St}\subseteq Pr^{L}$ we will denote the full $\infty$-subcategory generated by the stable and presentable $\infty$-categories 
\end{notation} 

The $\infty$-subcategory $Pr^{L}\subseteq\Cat$ is stable under colimits.
By $(iii)$ we can write the $\infty$-functor $Sp(-)$ as a colimit, so we deduce that the adjunction (\ref{eqLocLexEx}) continues to be valid if we restrict both to the $\infty$-categories to respective presentable case 
\begin{equation}
\label{eqPRLocLexEx}
\begin{tikzcd}
	{P_{r}^{L}} && {P_{r}^{L,St};}
	\arrow[""{name=0, anchor=center, inner sep=0}, "{Sp(-)}", shift left=3, from=1-1, to=1-3]
	\arrow[""{name=1, anchor=center, inner sep=0}, "i", shift left=2, hook', from=1-3, to=1-1]
	\arrow["\dashv"{anchor=center, rotate=-90}, draw=none, from=0, to=1]
\end{tikzcd}
\end{equation}
Using the symmetric monoidal structure of $Pr^{L}$, it is possible describe the $\infty$-functor \[Sp:Pr^{L}\to Pr^{L,St}\] as a tensor:   
\begin{equation}
\label{eqStabilizatorTensor}
\begin{split}
\C\otimes\mathcal{S}p\simeq & RFun(\C^{op},\mathcal{S}p)\simeq holim\{RFun(\C^{op},\mathcal{S}_*)\}\\ \simeq & holim \{\C\otimes\mathcal{S}_*\}\simeq holim \{\C_* \}\\ \simeq & \mathcal{S}p(\C).
\end{split}
\end{equation}

Now we explain how the monoidal structure of $Pr^{L}$ induces a monoidal structure on $Pr^{L,St}$.  

The $\infty$-category $\mathrm{Mod}_{Sp}(Pr^{L})$ has the symmetric monoidal structure of left $\Sp$-modules object of $Pr^{L}$ (i.e. the tensor is relative tensor product \cite[\S 4.4]{HA}).

It is known that there is a free-forgetful adjunction 
 \begin{equation} 
 \label{eqConfusione}
 -\otimes^{L}Sp:Pr^{L}\rightleftarrows \mathrm{Mod}_{Sp}(Pr^{L}):U.
  \end{equation}
and, \cite[Proposition 4.8.2.18]{HA}, Lurie proves that the forgetful $\infty$-functor $U:\mathrm{Mod}_{Sp}(Pr^{L})\to Pr^{L} $ is fully faithful and that its essential image is the $\infty$-category of presentable stable $\infty$-categories $Pr^{L,St}$. So the adjoint restricts to an adjoint equivalence
 \begin{equation}
 \label{eqPrstModSpPrL}
    -\otimes^{L}Sp:Pr^{L,St}\rightleftarrows \mathrm{Mod}_{Sp}(Pr^{L}):U.
\end{equation}
 
 We conclude that $Pr^{L,St}$ inherits the symmetric monoidal structure from the $\infty$-category of left $\Sp$-module objects of $Pr^{L}$, $\mathrm{Mod}_{Sp}(Pr^{L})$.
 We denote by $\otimes^{L,St}$ its tensor: 
\begin{equation}
\label{eqTensorPrlsTUstrega}
\begin{split}    
-\otimes^{L,St}-:Pr^{L,St} & \times Pr^{L,St}\to Pr^{L,St} \\ (\A & ,\B)\mapsto \A\otimes^{L,St}\B:=: U(-\otimes^{L}\Sp(\A)\; _{\Sp}\otimes_{\Sp} -\otimes^{L}\Sp(\B)) \\& \simeq Bar_{\Sp}(\A\otimes^{L}\Sp, \B\otimes^{L}\Sp). 
\end{split}
\end{equation}
and its unit is $\Sp$. 

Moreover, it preserves small colimits in each component, \cite[Corollary 4.4.2.15.]{HA}.
In short, $Pr^{L,St}$ is a commutative algebra object in $\Cat(\K)$, \cite[\S 3]{DoniHigherCategorical}, where $\K$ is the collection of all small simplicial sets.

As a trivial consequence, we obtain that the forgetful-free adjunction 
\[-\otimes^{L}\Sp:Pr^{L}\rightleftarrows \mathrm{Mod}_{\Sp}(Pr^{L}):U \] is monoidal.

\begin{warning}
 The $\infty$-functor $-\otimes^{L} \Sp$ in (\ref{eqConfusione}) is not the equal to $\Sp(-)$ in \eqref{eqStabilizatorTensor}.   
\end{warning}

\begin{notation}
Let $Pr^{L,St}_{cg}$ be the full $\infty$-subcategory spanned of $Pr^{L,St}$ spanned by stable and compactly generated $\infty$-categories.     
\end{notation}

The $\infty$-category inherits a symmetric monoidal structure via the inclusion $Pr^{L,St}_{cg}$ since the tensor \eqref{eqTensorPrlsTUstrega} sends compactly generated $\infty$-categories in compactly generated $\infty$-category. We denote this tensor as the tensor in \eqref{eqTensorPrlsTUstrega}.


\begin{example}
Now we want to prove that $Pr^{L,St}$ and $Pr^{L,St}_{cg}$ are left $Alg(\Sp)$-module objects of $\Cat(\K)$ and we need the following results (see \cite[Corollary 5.51]{DoniHigherCategorical} for the notation). 

\begin{lemma}[{\cite[Corollary 4.2.3.7]{HA}}]
\label{corFactorization}
The canonical projection 
\[ p_0: \mathrm{LMod}(\Sp)\to Alg(\Sp) \]
is both a cartesian and cocartesian fibration.

Moreover, both its cocartesian Grothendieck construction 
\[\begin{tikzcd}[ampersand replacement=\&]
    {\int p:} \& {Alg^{op}} \&\& \Cat \\
    \& E \&\& {\mathrm{LMod}_{E}} \\
    \& A \&\& {\mathrm{LMod}_{A}.}
    \arrow["f", from=2-2, to=3-2]
    \arrow["{f_{!}}", from=2-4, to=3-4]
    \arrow[from=1-2, to=1-4]
    \arrow[shorten <=16pt, shorten >=11pt, maps to, from=2-2, to=2-4]
\end{tikzcd}\]

and its cartesian Grothendieck construction 
\[\begin{tikzcd}[ampersand replacement=\&]
    {\int^{op} p:} \& {Alg^{op}} \&\& \Cat \\
    \& E \&\& {\mathrm{LMod}_{E}} \\
    \& A \&\& {\mathrm{LMod}_{A}}
    \arrow["f", from=2-2, to=3-2]
    \arrow["{f^*}"', from=3-4, to=2-4]
    \arrow[from=1-2, to=1-4]
    \arrow[shorten <=16pt, shorten >=11pt, maps to, from=2-2, to=2-4]
\end{tikzcd}\]
factor through the inclusion $\iota:Pr^{L,St}\to\Cat.$ 
\end{lemma}
\begin{proof}
As a consequence of \cite[Corollary 5.51 and Remark 5.50]{DoniHigherCategorical}, the $\infty$-functor $p_0$ is a cocartesian and cartesian fibration.

Moreover, its Grothendieck construction factors through $\iota$ because $\mathrm{LMod}_{E}$ is a stable presentable $\infty$-category for each $\mathbb{E}_1$-ring $E$, and the $\infty$-functors $f_{!}$ and $f^{*}$ preserve small colimits because both of them are left adjoint $\infty$-functors, \cite[Corollario 5.51]{DoniHigherCategorical}.   
\end{proof}

We denote by $\mathrm{LMod}_{(-)}$ the $\infty$-functor $\int p_0$ in the above lemma:  

\[\mathrm{LMod}_{(-)} : Alg\to Pr^{L,St}: \A\mapsto \mathrm{LMod}_{\A}.\]

This $\infty$-functor is monoidal, \cite[Corollary 4.8.5.13., Theorem 4.8.5.16. (4)]{HA}, and this implies that we can restrict the left action of $Pr^{L,St}$ on itself, and we obtain that $Pr^{L,St}$ is a left $Alg(\Sp)$-module object in $\Cat(\K)$ with the following tensor product: 

\begin{equation}
\label{eqTensorPrL}
-^{\Sp}\otimes^{L,St}-:Alg \times Pr^{L,St}\to Pr^{L,St}: (R,\A)\mapsto \mathrm{LMod}_{R}\otimes^{L,St} \A. 
\end{equation}

Moreover, since for each $\A$ in $Alg(\Sp)$ the $\infty$-category $\mathrm{LMod}_{\A}$ is compactly generated and stable, the above $\infty$-functor can be restricted further, and it defines a structure of left $Alg(\Sp)$-module objects of $\Cat(\K)$ on $Pr_{cg}^{L,St}$

\begin{equation}
\label{eqTensorPrLStcg}
-^{\Sp}\otimes^{L,St}-:Alg(\Sp) \times Pr_{cg}^{L,St}\to Pr_{cg}^{L,St}: (R,\A)\mapsto \mathrm{Mod}_{R}\otimes^{L,St} \A. 
\end{equation}

\end{example}



\section{$\Cat^{perf}$ and $P^{L,St}_{\omega}$: non-geometrical Morita theory  }
\label{subsecCatperfPRomega}

In this section, we recall and sometimes reproof some results in the pioneering work \cite{blumberg2013universal} which will be important in this paper. 

First, we define some $\infty$-categories.
In \Cref{exIdempo}, we defined the symmetric monoidal $\infty$-category $\Cat^{Idem}$ whose objects are (not necessarily small) idempotent complete $\infty$-categories and whose morphisms are $\infty$-functors (with no other properties),  
and in \Cref{exCatRIdem}, we defined the symmetric monoidal $\infty$-category $\Cat^{RIdem}$ whose objects are idempotent complete $\infty$-categories with finite colimits and whose morphisms are finite-colimits-preserving $\infty$-functor.

\begin{definition}
    We denote by $\Cat^{perf}\subseteq \Cat^{RIdem}$ the full $\infty$-subcategory whose objects are small stable idempotent complete $\infty$-categories.
\end{definition}

\begin{definition}[{\cite[Definition 5.4.2.19]{HTT}}]
We denote by $Acc_{\omega}\subseteq\widehat{\Cat}$ the $\infty$-subcategory defined as follows:
\begin{itemize}
    \item the objects are $\omega$-accessible $\infty$-categories;

    \item an $\infty$-functor $F:\C\to\D$ between accessible $\infty$-categories belongs to $Acc_{\omega}$ if and only if it preserves $\omega$-filtered colimits and preserves compact objects.
    
\end{itemize}
\end{definition}

\begin{notation}
    We denote by $Pr^{L}_{\omega}$ the $\infty$-subcategory of $Acc_{\omega}$ spanned by compactly generated $\infty$-categories and whose morphisms are preserving-small-colimits $\infty$-functors. 
\end{notation}

In \cite[Proposition 5.4.2.17]{HA}, the author proves that there exists an equivalence
\begin{equation}
    \label{eqZZoCA}
 Ind_{\omega}(-): \Cat^{Idem}\rightleftarrows Acc_{\omega}:(-)^{\omega}. 
 \end{equation}

By \Cref{rmkCatIdemsymmetricMonoidal}, we know that $\Cat^{Idem}$ has a symmetric monoidal structure, so $Acc_{\omega}$ inherits a symmetric monoidal structure.

By \Cref{rmkCatIdemsymmetricMonoidal}, there is an inclusion $\iota^{RIdem}:\Cat^{RIdem}\to \Cat^{Idem}$, and we restrict \eqref{eqZZoCA} along this inclusion. Using \cite[Theorem 5.5.1 (3)]{HTT}, we obtain that the essential image of $\Cat^{RIdem}$ via $Ind_{\omega}(-)$ is $Pr^{L}_{\omega}$. So we obtain the next equivalence:
\begin{equation}
\label{eqGaFi}
Ind_{\omega}(-): \Cat^{RIdem}\rightleftarrows Pr^{L}_{\omega}:(-)^{\omega}. 
\end{equation}

By \Cref{rmkCatIdemsymmetricMonoidal}, we know that $\Cat^{RIdem}$ has a symmetric monoidal structure, so $Pr^{L}_{\omega}$ inherits a symmetric monoidal structure. For the moment we denote the tensor in $Pr^{L}_{\omega}$ by $-\otimes^{L,\omega}-$:
\[-\otimes^{L,\omega}-:(\A,\B)\mapsto  \A\otimes^{L,\omega}\B := Ind(\A^{\omega}\otimes^{RIdem} \B^{\omega}).\] 

Now we prove that the inclusion $Pr^{L}_{\omega}\to Pr^{L}$ is monoidal.

Let $\F$ be the collection of finite simplicial sets, let $N$ be the collection as in \Cref{exIdempo}, and $\K$ the collection of all small simplicial sets. From the last part of \Cref{secAddingColimits}, we have the following (a priori non-commutative) diagram:

\[\begin{tikzcd}[ampersand replacement=\&]
	{\Cat^{rex}} \& {\Cat^{RIdem}} \& {Pr^{L}_{\omega}} \\
	\&\& {Pr^{L}} \\
	\&\& {\Cat(\K).}
	\arrow[""{name=0, anchor=center, inner sep=0}, "{P^{N}_{\F}}", shift left=2, from=1-1, to=1-2]
	\arrow[""{name=1, anchor=center, inner sep=0}, "Ind", shift left=2, from=1-2, to=1-3]
	\arrow[hook', from=1-3, to=2-3]
	\arrow["{P^{\K}_{\F}(-)}"', from=1-1, to=2-3]
	\arrow[""{name=2, anchor=center, inner sep=0}, "{(-)^{\omega}}", shift left=2, from=1-3, to=1-2]
	\arrow[""{name=3, anchor=center, inner sep=0}, shift left=2, hook', from=1-2, to=1-1]
	\arrow[hook', from=2-3, to=3-3]
	\arrow["{P^{\K}_{\F}(-)}"', curve={height=18pt}, from=1-1, to=3-3]
	\arrow["\cong"{description}, draw=none, from=1, to=2]
	\arrow["\dashv"{anchor=center, rotate=-90}, draw=none, from=0, to=3]
\end{tikzcd}\]

To achieve our target, we need to prove the following result: for each object $\C\in \Cat^{rex}$, there exists an equivalence of $\infty$-categories 
\begin{equation}
\label{eqfinalmenteCazzo}
    P^{\K}_{\F}(\C)\simeq Ind(P^{N}_{\F}(\C)).
\end{equation}
We prove the result by showing that $ Ind(P^{N}_{\F}(\C))$ satisfies the properties of $P^{\K}_{\F}(\C)$, \Cref{defPKH}.
Property $(1)$ is trivial consequence of \Cref{exTHeTrick}.
Let $\B\in \Cat(\K)$ be an $\infty$-category which admits all small colimits; this implies, in particular, that it is idempotent complete, see \Cref{exIdempo}.

The following chain of equivalences holds: 
\begin{equation}
\begin{split}
Fun_{\Cat(\K)}(Ind (P^{N}_{\F}(\C)),B) \simeq & Fun_{\Cat^{RIdem}}((Ind(P^{N}_{\F}(\C)))^{\omega},B) \\ \simeq & Fun_{\Cat^{RIdem}}(P^{N}_{\F}(\C),B)\\ \simeq & Fun_{\Cat^{rex}}(P^{N}_{\F}(\C),B).
\end{split}
\end{equation}

The first equivalence holds by \cite[Proposition 5.5.1.9]{HTT}, the second one because of \eqref{eqGaFi}, and the last one because $P^{N}_{\F}(-): \Cat^{rex}\to \Cat^{RIdem}$ is a localization, \Cref{rmkLocalizationrexRIdem}.

Now we can prove our target.
The following chain of equivalences holds:

\begin{equation}
    \begin{split}
        \A\otimes^{L,\omega}\B  = & Ind(\A^{\omega}\otimes^{RIdem} \B^{\omega}) \cong Ind(P^{N}_{\F}(\A^{\omega})\otimes^{RIdem} P^{N}_{\F}(\B^{\omega})) \\
        \cong & Ind(P^{N}_{\F}(\A^{\omega}\otimes^{rex}\B^{\omega})) \overset{\dagger}{\cong} P_{\F}^{\K}(\A^{\omega}\otimes^{rex}\B^{\omega})) \\ \cong & P_{\F}^{\K}(\A^{\omega})\otimes^L P_{\F}^{\K}(\B^{\omega}) \overset{\square}{\cong} Ind(\A^{\omega})\otimes^L Ind(\B^{\omega}) \\ \cong &  \A\otimes^{L}\B.  
    \end{split}
\end{equation}
The above equivalences are valid because $P^{\mathscr{K}}_{\F}(-)$ and $P^{N}_{\F}(-)$ are monoidal, 
$P^{K}_{\F}(-):\Cat^{rex}\rightleftarrows \Cat^{RIdem}$ is a localization by \Cref{rmkLocalizationrexRIdem}, the $\square$ holds by \Cref{exTHeTrick},
the $\dagger$ follows by the equivalence \eqref{eqfinalmenteCazzo}, 
and moreover, we use that the monoidal structure of $Pr^{L}$ is the restriction of the monoidal structure of $\Cat(\K)$, see \cite[Proposition 4.8.1.15]{HA}.

We conclude that the (non-full) inclusion $Pr_{\omega}^{L}\to Pr^{L}$ is monoidal and that the tensor product and the unit of $Pr_{\omega}^{L}$ coincide with the unit and the tensor of $Pr^{L}$. For this reason, from now on, we will denote by $\otimes^{L}$ the tensor of $Pr_{\omega}^{L}$, which was previously denoted by $\otimes^{L,\omega}$.
So the triple $(Pr^{L}_{\omega},\otimes^{L},\Top)$ forms a monoidal structure: 
\begin{equation}
\label{eqTensorPrlOmega}
\begin{split}    
-\otimes^{L}-:Pr^{L}_{\omega} & \times Pr^{L}_{\omega}\to Pr^{L}_{\omega}\\ (\A & ,\B)\mapsto \A\otimes^{L}\B:=RFun(\A^{op},\B). 
\end{split}
\end{equation}

In \cite[Theorem 5.5.3.18, Proposition 5.5.7.6]{HTT}, the author proves that $Pr^{L}_{\omega}$ admits small colimits and that the small colimits in $Pr^{L}_{\omega}$ are computed in $Pr^{L}$. Since the tensor product $\otimes^{L}$ of $Pr^{L}$ preserves small colimits, the tensor product $\otimes^{L}$ of $Pr^{L}_{\omega}$ (which is the same as $Pr^{L}$) also preserves them.
In short, $Pr^L_{\omega}$ is a commutative algebra object in $\Cat(\K)$, where $\K$ is the collection of all small simplicial sets.

Now, we restrict the equivalence in \eqref{eqPrstModSpPrL} via the inclusion $Pr_{\omega}^{L,St}\to Pr^{L,St} $ and we obtain the equivalence
\[Pr^{L,St}_{\omega}\rightleftarrows \mathrm{Mod}_{\Sp}(Pr^{L}_{\omega}),\]
since $\Sp$ is compactly generated \cite[Proposition 1.4.3.7]{HA}.
The $\infty$-category $\mathrm{Mod}_{Sp}(Pr_{\omega}^{L})$ has the symmetric monoidal structure of left $\Sp$-module objects of $Pr^{L}_{\omega}$ (i.e. the tensor is the relative tensor product, \cite[\S 4.4]{HA}).

 We conclude that $Pr^{L,St}_{\omega}$ inherits the symmetric monoidal structure from the $\infty$-category of left $\Sp$-module objects of $Pr^{L}_{\omega}$, $\mathrm{Mod}_{Sp}(Pr^{L}_{\omega})$. Reasoning as before, we obtain that the (non-full) inclusion
\[Pr^{L,St}_{\omega}\to Pr^{L,St}\] is monoidal and then the tensor and the unit of $Pr^{L,St}_{\omega}$ are the same as those of $Pr^{L,St}$. Therefore, we use the same notation to denote them: 
\begin{equation}
\label{eqTensorPrlsT}
\begin{split}    
 -\otimes^{L,St}-:Pr^{L,St}_{\omega} & \times Pr_{\omega}^{L,St}\to Pr_{\omega}^{L,St} \\ (\A & ,\B)\mapsto \A\otimes^{L,St}\B:= U(-\otimes^{L}\Sp(\A)\; _{\Sp}\otimes_{\Sp} -\otimes^{L}\Sp(\B)) \\ & \simeq Bar_{\Sp}(\A\otimes^{L}\Sp, \B\otimes^{L}\Sp).
\end{split}
\end{equation}

and its unit is $\Sp$. Moreover, the above tensor preserves small colimits in each component, \cite[Corollary 4.4.2.15.]{HA}.
In short, $Pr^{L,St}_{\omega}$ is a commutative algebra of $\Cat(\K)$ where $\K$ is the collection of all small simplicial sets.

Finally, the inclusion $\iota: \Cat^{perf}\to \Cat^{RIdem}$ is fully faithful (in both the $\infty$-categories the morphism are preserving-finite-colimits $\infty$-functor), then we restrict the equivalence $(\ref{eqGaFi})$ via $\iota$ and we find the adjoint equivalence: 
\begin{equation}
\label{eqCatPerfPrStLOmega}
Ind(-):\Cat^{perf}\rightleftarrows Pr^{L,St}_{\omega}:(-)^{\omega}.
\end{equation}
This above equivalence was previously described in \cite{blumberg2013universal}.

\begin{remark}[{Technical point}]
\label{rmktecnicalPoint}
Let $\C$ be an object of $Pr^{L,St}_{\omega}$. The full $\infty$-subcategory generated by compact objects $C^{\omega}\subseteq C$ is not necessarily a small $\infty$-category, (\cite[Remark 5.4.2.11]{HA}), but it is only essentially small. So in \Cref{eqCatPerfPrStLOmega}, the $\infty$-functor $(-)^{\omega}$ actually carries an object $\C\in Pr^{L,St}_{\omega}$ into a small $\infty$-category equivalent to $\C^{\omega}$.
\end{remark}

We conclude that $\Cat^{perf}$ inherits a symmetric monoidal structure from $Pr^{L}_{\omega}$.
 
\begin{equation}
\label{eqTensorPerf}
\begin{split}    
-\otimes^{perf}-:\Cat^{perf} & \times \Cat^{perf}\to \Cat^{perf} \\ (\A & ,\B)\mapsto (Ind(\A)\otimes^{L,St}Ind(\B))^{\omega}.
\end{split}
\end{equation}
and its unit is $\Sp^{\omega}$. Moreover, it preserves small colimits in each component because the tensor of $Pr_{\omega}^{L,St}$ does so.
In short, $\Cat^{perf}$ is a commutative algebra object of $\Cat(\K)$ where $\K$ is the collection of all small simplicial sets.

\begin{remark}
In \cite{Cohn}, the proof that $\Cat^{perf}$ has a monoidal structure is wrong. The author uses the method in \cite[Remark 2.2.1.2]{HA}
on the inclusion $Pr_{\omega}^{L,St}\to Pr_{\omega}^{L,St}$. However, this criterion cannot be applied because this inclusion is non-full.

In \cite{blumberg2013universal}, the authors construct the monoidal structure, but they do not write a complete proof.
\end{remark}


In \cite[Corollary 4.25]{blumberg2013universal}, the authors prove that $\Cat^{perf}$ and $Pr^{L,St}_{\omega}$ are compactly generated $\infty$-categories, and then they are presentable $\infty$-categories. So, they are actually commutative algebra objects of $Pr^{L}$ (actually, of $Pr^{L}_{cg}$). In particular, this implies that both have an internal-hom.

For us, it is important to find that $\Cat^{perf}$ and $Pr^{L,St}_{\omega}$ are left $Alg$-module objects of $Pr^{L}$.

\begin{proposition}
\label{propequivalenceCatperfprlomegaalgspmod}
The $\infty$-categories $\Cat^{perf}$ and $Pr^{L,St}_{\omega}$ are two left $Alg$-module objects of $Pr^{L}$.
Moreover, $\Cat^{perf}$ and $Pr^{L,St}_{\omega}$ are equivalent as left $Alg$-module objects of $Pr^{L}$.
\end{proposition}
\begin{proof}
Since the monoidal structure in $\Cat^{perf}$ and $Pr^{L,St}_{\omega}$ is induced by the monoidal structure in $Pr_{cg}^{L,St}$, we can restrict the action of $Alg$ on $Pr_{cg}^{L,St}$ given in \eqref{eqTensorPrLStcg} to the $\infty$-categories $\Cat^{perf}$ and $Pr^{L,St}_{\omega}$. Since there is a monoidal equivalence $\Cat^{perf}\simeq Pr^{L,St}_{\omega}$, they are equivalent as left $Alg$-module objects of $Pr^{L}$.
\end{proof}

In the same way, using the monoidal equivalence $\Cat^{perf}\simeq Pr^{L,St}_{\omega}$, it is possible to see $\Cat^{perf}$ as a left $Pr^{L,St}_{\omega}$-module of $Pr^{L}$ and, vice versa, $Pr_{\omega}^{L,St}$ as a left $\Cat^{perf}$-module of $Pr^{L}$.

\begin{remark}
\label{rmkSptensorizationonCatPerfPrSt}
Now we sum up the higher structures of $\Cat^{perf}$, $Pr^{L,St}_{\omega}$ and we fix the notation for all of them:
\begin{itemize}
    \item $\Cat^{perf}$ as left $Alg$-module objects of $Pr^{L}$:
    \[-^{\Sp}\otimes^{perf}-:Alg \times \Cat^{perf}\to \Cat^{perf}:(m,v)\mapsto (LMod_{m}\otimes^{L,St}Ind(v))^{\omega};  \]

    \item $\Cat^{perf}$ as associative algebra of $Pr^{L}$:
    \[-\otimes^{perf}-:\Cat^{perf} \times \Cat^{perf}\to \Cat^{perf}:(m,v)\mapsto (Ind(m)\otimes^{L,St}Ind(v))^{\omega};  \]

    \item $\Cat^{perf}$ as left $Pr_{\omega}^{L,St}$-module objects of $Pr^{L}$:
    \[-^{\omega}\otimes^{perf}-: Pr_{\omega}^{L,St} \times \Cat^{perf} \to \Cat^{perf}:(m,v)\mapsto (m\otimes^{L,St}Ind(v))^{\omega}.  \]

    \item $Pr_{\omega}^{L,St}$ as left $Alg$-module of $Pr^{L}$:
    \[-^{\Sp}\otimes^{L,St}-:Alg \times Pr^{L,St}_{\omega}\to Pr^{L,St}_{\omega}:(m,v)\mapsto LMod_{m}\otimes^{L,St}v;  \]

    \item $Pr_{\omega}^{L,St}$ as left $\Cat^{perf}$-module of $Pr^{L}$:
    \[-^{perf}\otimes^{L,St}-:\Cat^{perf} \times Pr^{L,St}_{\omega}\to Pr^{L,St}:(m,v)\mapsto Ind(m)\otimes^{L,St}v.  \]

    \item $Pr_{\omega}^{L,St}$ as associative algebra of $Pr^{L}$:
    \[-\otimes^{L,St}-:Pr_{\omega}^{L,St} \times Pr_{\omega}^{L,St} \to Pr_{\omega}^{L,St}:(m,v)\mapsto m\otimes^{L,St}v.  \]
\end{itemize}    
\end{remark}


\section{The $\infty$-categories of dg-categories}
\label{secDg}

In this section, we will explain the theory of dg-categories using enriched $\infty$-categories theory. This approach is not classical. Historically, dg-categories have been studied using the enriched model category theory. However, since Haugseng proved that enriched $\infty$-category theory is a flexible version of enriched model category theory, this modern approach can be useful and natural, see \cite[Theorem 5.8]{HaugsengRect} for the original result or \cite[\S 3.3]{DoniPhDThesis}.

Nowadays, the theory of enriched $\infty$-category is not well-developed and we must make some choices. We did them in \cite[\S 5]{DoniHigherCategorical} where we explain also our convention. We suggest reading \cite[\S 5]{DoniHigherCategorical} before approaching the remainder of the article. 

This section should, however, be seen as \enquote{the algebraic geometric way to see a $k$-linearization}, which in practice means that the $\infty$-categories are at least $\Z$-linear from the start. This has the drawback that when we add properties to an $\infty$-category, we need to consider the enrichment. For example, adding small colimits in this context is not as easy as for an $\infty$-category (see \cite[\S 5.1.5]{HTT}), because adding colimits means adding enriched colimits. Nowadays, as far as we know, it is not possible to add freely enriched colimits.

In the first \Cref{subsecDgcatLmod}, we will prove that $\D(k)$-enriched $\infty$-categories are $\LModk$-enriched $\infty$-categories. This result was proved by Lurie in \cite{HA} and by Haugseng-Gepner in \cite{GepHauEnriched}, and it was precedently proved by using  enriched model category theory, \cite[\S 3]{DoniPhDThesis}.

In \cite{DoniHigherCategorical}, we proved that $\LModk$-enriched $\infty$-categories are left $\Hk$-module objects of $\Cat^{\Sp}$. This result is very important for this work because it allows us to redefine a $k$-linearization in, at least, three ways:
\begin{itemize}
    \item an $\LModk$-enriched $\infty$-category. This is the algebro geometric view;
    \item or, a $ALg$-enriched $\infty$-functor from the $Alg$-enriched $\infty$-category $\un{\Hk}$ to the $Alg(\Sp)$-left tensored $\infty$-category $\Cat^{\Sp}$. This new view is the starting point for a general notion of $k$-linearization, see \Cref{subsecGeneraleklinear} and it is the spectral algebraic view;
    \item or, a left $\Hk$-module object of $\Cat^{\Sp}$.
\end{itemize}

In \Cref{subsecpretDGCat}, we will develop the Morita theory for $\D(k)$-enriched $\infty$-categories.


\subsection{From $\Cat^{\D(k)}$ to $\Cat^{\LModk}$}
\label{subsecDgcatLmod}

In this subsection, we will explain how to define an equivalence:
\begin{equation}
\label{eqDgSimeqEnModk}
\Cat^{\D(k)}\simeq\Cat^{\LModk},
\end{equation}
which will be a change of enrichment (also known as a change of base).

To achieve our result, we apply the construction in \cite[Notation 5.19]{DoniHigherCategorical} to a symmetric monoidal equivalence
\begin{equation}
\label{eqChangeEnr}
\hat{F}:\D(k)\simeq \LModk.
\end{equation}
In Theorem \cite[Theorem 7.1.2.13]{HA}, the author proves the existence of this symmetric monoidal equivalence. For the sake of completeness, we recall how the $\infty$-functor $\hat{F}$ is defined.

This $\infty$-functor is constructed in three steps:

First, we construct a right $t$-exact $\infty$-functor
\begin{equation}
\label{eqDmenok}
    F:\D^-(k)\to \LModk.
\end{equation}

The $\infty$-category $\LModk$ has a right complete $t$-structure, \Cref{exSperiamoUltimoDaCitare}, then we can use the universal property of the bounded above derived $\infty$-category $\D^-(k)$, \cite[Theorem 1.3.3.2]{HA}, which states that the following equivalence holds:
\begin{equation}
 \label{DmenoUniProp}   
 Fun^{p,\heartsuit}(\D^-(k), \LModk) \simeq N(\Hom(k\text{-}\mathrm{Mod},k\text{-}\mathrm{Mod})).
 \end{equation} 
 Here, $Fun^{p,\heartsuit}(\D^-(k), \LModk)$ is the $\infty$-category of right $t$-exact $\infty$-functors that carry projective objects of $\D^-(k)^{\heartsuit}\simeq k\text{-}\mathrm{Mod}$ into the heart of $\LModk$.

\begin{remark}
\Cref{DmenoUniProp} is truly interesting. It implies that any ordinary derived functor can be lifted to a right $t$-complete $\infty$-functor between derived $\infty$-categories. Indeed, this equivalence carries an object $F\in Fun^{p,\heartsuit}(\D^-(k), \LModk)$ into the derived functor \[k\text{-}\mathrm{Mod}\xrightarrow{\deg_{0}}Ch_{\geq 0}(k)\xrightarrow{ \mathbb{L}hF }Ch_{\geq 0}(k)\xrightarrow{H_{0}} k\text{-}\mathrm{Mod} \] of a right exact functor $hL$. See \cite[\textsection 2.3]{riehl2014categorical} for an introduction to derived functors.
\end{remark}

In our case, we lift the identity functor $id_{k\text{-}\mathrm{Mod}}$, and we obtain the desired $\infty$-functor \eqref{eqDmenok}. This concludes the first step.

The Proposition \cite[Proposition 1.3.5.24]{HA} allows us to identify $\D^{-}(k)$ with the full $\infty$-subcategory of $\bigcup_{n\geq 0}\D(k)_{\geq -{n}}\subseteq \D(k)$ spanned by the right bounded objects. Since $F$ is a right $t$-exact $\infty$-functor, $F$ sends right bounded objects into right bounded objects. This implies that we can restrict it to every right bounded component and define a diagram $F_{*}:\D(k)_{\geq *}\to\LModk^{\geq *}$ as in the figure below:

\[\begin{tikzcd}
	{} && {} \\
	{\D(k)_{\geq -2}} && {\LModk^{\geq -2}} \\
	{\D(k)_{\geq -1}} && {\LModk^{\geq -1}} \\
	{\D(k)_{\geq 0}} && {\LModk^{\geq 0}.}
	\arrow["{\tau_{\geq 0}}", from=3-1, to=4-1]
	\arrow["{F_0}", from=4-1, to=4-3]
	\arrow["{F_{-1}}", from=3-1, to=3-3]
	\arrow["{\tau_{\geq 0}}", from=3-3, to=4-3]
	\arrow["{\tau_{\geq -1}}", from=2-3, to=3-3]
	\arrow[dashed, no head, from=2-3, to=1-3]
	\arrow["{\tau_{\geq -1}}", from=2-1, to=3-1]
	\arrow[dashed, no head, from=2-1, to=1-1]
	\arrow["{F_{-2}}", from=2-1, to=2-3]
\end{tikzcd}\]

We take the limit of $F_{*}$, and since $\D(k)$ and $\LModk$ have right complete $t$-structure, we obtain our desired $\infty$-functor \eqref{eqChangeEnr}:
\begin{equation}
\label{eqD(k)LModK}
\text{lim}(\D(k)_{\geq *})\simeq \D(k)\xrightarrow{ \hat{F}= \text{lim} F_{*}} \text{lim}(\LModk^{\geq *})\simeq \LModk.
\end{equation}
That concludes the second step.
In the third step, we need to prove that $\hat{F}$ is monoidal and equivalence. This is achieved using Shipley-Schwede's recognition principle \cite[\textsection 7.1.2]{HA}, keeping in mind the following two facts: 
\begin{itemize}
    \item the unit object of the symmetric monoidal structure of $\D(k)$ is \[k[0]:=\dots\to 0\to k \to 0\to \dots \] and $k[0]$ is compact and corepresents the homology (in \cite[Theorem 7.1.2.1]{HA}, Lurie calls an object as $k[0]$ a generator of $D(k)$) of chain complexes:
    \begin{equation}
    \begin{split}
        \pi_{0}\Hom_{\D(k)}(k[0],[-n](-))\simeq & \pi_{0}\Hom_{\D(k)}(k[n],-) \\ & \simeq Ext_{\D(k)}^{-n}(k[0],-)  \simeq H_{n}(-):h\D(k)\to Set;
    \end{split}
    \end{equation}

    \item $\Hk$ is the endomorphism algebra of $k[0]$. Now, we prove this. Since $\D(k)$ is a presentable stable $\infty$-category, it can be considered as an $\Sp$-enriched $\infty$-category (see \cite[\S 7]{GepHauEnriched} or \cite{heine2023equivalence}). It is easy to see that the endomorphism object of $k[0]$ is the spectrum object $\D(k)(k[0],k[0])$, whose representation as a spectrum is: 
    \[
    (\dots,  \Hom_{\D(k)}((k[1],k[0]),\Hom_{D(k)}(k[0],k[0]),\Hom_{\D(k)}(k[-1],k[0])   ,\dots),
    \]
    for a proof of this easy fact, see \cite[Theorem 4.8.5.16]{HA} and \cite[\S 4.8.5]{HA}.
    $\D(k)(k[0],k[0])$ is a spectrum object because, by the universal property of pullbacks, the following equivalence holds: \[\Omega\Hom_{\D(k)}(k[n],k[0])\simeq\Hom_{\D(k)}(\Sigma k[n],k[0])\simeq \Hom_{\D(k)}(k[n+1],k[0])\].
    In particular, the homotopy groups of this spectrum are 
    \begin{equation}
    \begin{cases}
        \pi_{0}\D(k)(k[0],k[0])= \pi_{0}\Hom_{\D(k)}(k[0],k[0])\simeq H_{0}k[0]\simeq k & \text{if } n= 0; \\
        \pi_{n}\D(k)(k[0],k[0])= \pi_{0}\Hom_{\D(k)}(k[n],k[0])\simeq H_{n}k[0]\simeq 0 & \text{if} n\neq 0,
    \end{cases}
    \end{equation}
    We conclude that $\D(k)(k[0],k[0])$ is equivalent to the Eilenberg-Maclane spectrum $\Hk$. 
\end{itemize}

Finally, using the functoriality of the $\infty$-functor $\Cat^{(-)}$, \cite[Equation (37)]{DoniHigherCategorical}, we obtain the desired \Cref{eqDgSimeqEnModk}.
So, above we proved the following result:
\begin{proposition}
\label{propEqDkenrLModkenric}
The $\infty$-category of $\D(k)$-enriched  $\infty$-categories $\Cat^{\D(k)}$ and the $\infty$-category of $\LModk$-enriched $\infty$-categories $\Cat^{\LModk}$ are equivalent. 
\end{proposition}


\subsection{Morita theory for dg-categories }
\label{subsecpretDGCat}

First of all, we recall what a dg-category over $k$ is for us.
\begin{definition}
A $dg$-category over $k$ is an object of the $\infty$-category of $\D(k)$-enriched $\infty$-categories, denoted as $\Cat^{\D(k)}$. 
\end{definition}

Sometimes, with abuse, we say dg-category instead of dg-category over $k$.

\begin{remark}
Our definition is equivalent to the classical definition, where a $dg$-category is a $Ch(k)$-enriched category, e.g. \cite[Definition 1.1]{canonaco2017tour}, since $\Cat^{\D(k)}$ is the underlying $\infty$-category of the Dwyer-Kan model category of $Cat(Ch(k))$, detailed in \cite[\S 3]{DoniPhDThesis}.    
\end{remark}

\begin{remark}
In \Cref{subsecDgcatLmod}, we proved that $\Cat^{\D(k)}\simeq \Cat^{\LModk}$, so a dg-category over $k$ is also an object of $\Cat^{\LModk}$.

From now on we call dg-categories over $k$ also the objects of $\Cat^{\LModk}$.

Hereafter, we will use the $\infty$-category $\Cat^{\LModk}$ instead of $\Cat^{\D(k)}$.

\end{remark}

In this subsection, we enrich the results in \Cref{subsecCatperfPRomega}.
We begin with the $\Sp$-enriched case. By the end of this thesis, it will be evident that the $\Sp$-enriched case and the non-geometrical case in \Cref{subsecCatperfPRomega} are equivalent. Thus, we sometimes refer to the $\Sp$-enriched case as the non-geometrical case.

After analyzing the non-geometrical case, we will address the $k$-linear one, i.e., the $\LModk$-enriched case.

\begin{notation}
We refer to objects of $\Cat^{\Sp}$ as \textit{spectral $\infty$-categories}.
\end{notation}

\begin{notation}[{\cite[Example 5.25]{DoniHigherCategorical}}]
\label{notunderlyingspectralcategory}
    Let $\A\in\Cat^{\LModk}$ be a dg-category, we call \textit{the underlying spectral $\infty$-category of $\A$} the spectral $\infty$-category $U_{!}(\A)$ where $U_{!}: \Cat^{\LModk}\to \Cat^{\Sp}$ is the canonical change of enrichment. We denote it by $\A_{sp}$. 
\end{notation}

\begin{notation}[{\cite[Example 5.22,Notation 5.23]{DoniHigherCategorical}}]
We define the \textit{underlying $\infty$-category $\infty$-functor} as the following $\infty$-functor:
\begin{equation}
 \label{eqHomotopyHigherFunctorCat}  (-)_{o}:\Cat^{\V}\xrightarrow{Hom(\I_{\V},-)_{!}}\Cat^{\Top}\simeq \Cat:\;
 \A\mapsto \A_{o}
\end{equation}
Let $\A$ be an enriched $\infty$-category we call $\A_o$ its underlying $\infty$-category.  
Moreover, we denote by $hA$ its homotopy category $h\A_{o}$ and we call it \textit{the homotopy category of $\A$}. 
Let $f:\A\to\B$ be an object in $\Hom_{\Cat^{\V}}(\A,\B)$, we call $ f_o:\A_{o}\to\B_{o}$ the \textit{underlying $\infty$-functor of $f$}.
\end{notation}

We now consider the analog of the stable condition for spectral $\infty$-categories, namely the pretriangulated condition. It will be immediately clear that they are analogous to each other.

\begin{definition}[{\cite[Definition 9.1]{heine2023equivalence}}]
\label{defPretriang}
Let $\A\in \Cat^{\Sp}$ be a spectral $\infty$-category. We say that it is \textit{pretriangulated} if:
\begin{enumerate}
    \item ($\mathbf{0}$-element) $\A_{o}$ has a zero object $\mathbf{0}$;
    \item (mapping cone) $\A_{o}$ is closed under cofibers;
    \item (shift) The suspension $\infty$-functor $\Sigma:\A_{o}\to \A_{o}$ is well-defined, where $\Sigma(X) = \text{colim}(\mathbf{0}\leftarrow X\to \mathbf{0})$, and $\Sigma$ is an equivalence of $\infty$-categories.
\end{enumerate} 
and , for each object $Y\in\A$, the $\infty$-functors $\A(Y,-):\A_{o}\to \Sp$ and $\A(-,Y):\A_{o}^{op}\to \Sp$ are exact.
\end{definition}

\begin{remark}
In classical literature, dg-categories satisfying the properties of \Cref{defPretriang} are called \textit{strong pretriangulated}. Those satisfying only the properties (mapping cone) and (shift) are termed pretriangulated. We opt to call pretriangulated those dg-categories with all three properties since with this choice mapping cone and shift are (homotopy) colimits in the same dg-categories, which is significant to us. Using the classical convention for pretriangulated dg-categories (not necessarily strong), mapping cone and shift are defined as enriched (homotopy) colimits in the category of dg-modules.
\end{remark}

With a trick of category theory, it is easy to prove that the adjoint inverse of $\Sigma$ is the loop $\infty$-functor\[\Omega:\A_{o}\to \A_{o}:Y\mapsto lim(\mathbf{0}\to Y\leftarrow \mathbf{0}).\]
In \cite[Remark 9.7]{heine2023equivalence}, the author proves that $\Sigma\dashv \Omega$ is actually an $\Sp$-enriched adjunction. This is a consequence of the condition that, for each object $Y\in\A$, the $\infty$-functors $\A(Y,-):\A_{o}\to \Sp$ and $\A(-,Y):\A_{o}^{op}\to \Sp$ are exact.

Using \Cref{defStable}, we can rewrite \Cref{defPretriang} briefly as follows:

\begin{definition}
 \label{defPretriang1}
Let $\A$ be a spectral $\infty$-category. We say that it is \textit{pretriangulated} if its underlying $\infty$-category $\A_{o}$ is stable, and for each object $Y\in\A_{o}$, the $\infty$-functors $\A(Y,-):\A_{o}\to \Sp$ and $\A(-,Y):\A_{o}^{op}\to \Sp$ are exact.
\end{definition}

\begin{remark}
\label{rmkEnrichedcolimitsEnrichedLimits}
Let $\A$ be a pretriangulated spectral $\infty$-category, so $\A_{o}$ is stable and in particular $\A_{o}$ has all finite colimits.
Since in \Cref{defPretriang} and \Cref{defPretriang1} $\A_{o}(-,X):\A^{op}\to \Sp$ is exact, it preserves finite limits. This implies that $\A$ has $\Sp$-enriched colimits.
Indeed, let $d:I\to \A_{o}$ be a diagram with source a finite simplicial sets $I$ and $colim_{i\in I}d\in \A_o$ be the nadir of the colimit, then, for each $X\in\A_{o}$, the following chain of equivalences holds:
\begin{equation}
\label{eqColimtsLimitsLimits}
\A(colim_{i\in I} d(i),X)\simeq \A(lim_{i\in I^{op}}d^{op}(i),X)\simeq lim_{i\in I^{op}}\A(d(i),X),
\end{equation}
where the first holds because a colimit of $d$ in $\A_o$ is a limit of $d^{op}$ in $\A_o^{op}$ (see \cite[Proposition 12.1.7]{RiehlElements}) and the second one because $\A(-,X)$ preserves finite limits.
\eqref{eqColimtsLimitsLimits} is the $\infty$-categorical analog of the universal property of the enriched (or, conical) colimits of the ordinary case, see \cite[\S A.5]{RiehlElements}. 
In similar way the condition that $\A(X,-):\A_{o}\to \Sp$ preserves limits implies that $\A$ admits $\Sp$-enriched limits.

From what we know, the theory of enriched (co)limits is not well defined.
\end{remark}

\begin{remark}
\label{remarkCAZZOCAZZOCAZZO}
From the last definition, it is clear that the stability condition depends almost entirely on the underlying $\infty$-category $\A_o$ of a $\Sp$-enriched $\infty$-category.
In \cite[\S 9]{heine2023equivalence}, the author proves that the stability condition for a $\Sp$-enriched $\infty$-category is completely contained in the underlying $\infty$-category $\A_o$.
Indeed, he defines the full $\infty$-subcategory $\Cat^{\Sp,St}\subseteq \Cat^{\Sp}$ spanned by pretriangulated spectral $\infty$-categories and, in \cite[Corollary 9.21]{blumberg2013universal}, he proves that the underlying $\infty$-category $\infty$-functor define an equivalence:
    \[(-)_{o}:\Cat^{\Sp,St}\simeq \Cat^{Ex}.\]
    Moreover, he constructs a localization
    \begin{equation}
    \begin{tikzcd}[ampersand replacement=\&]
	{\Cat^{\Sp}} \&\& {\Cat^{\Sp,St}}
	\arrow[""{name=0, anchor=center, inner sep=0}, shift left=3, hook', from=1-3, to=1-1]
	\arrow[""{name=1, anchor=center, inner sep=0}, "{\mathscr{S}}", shift left=3, from=1-1, to=1-3]
	\arrow["\dashv"{anchor=center, rotate=-90}, draw=none, from=1, to=0]
\end{tikzcd}
\end{equation}

whose right adjoint is lax monoidal. This result is interesting because it is the $\infty$-categorical analogous to the triangulated localization $\Psi_{Tri}(-)$ defined in \cite{blumberg2013universal}.
\end{remark}

\begin{remark}
From \Cref{defPretriang1}, we obtain some obvious consequences:
\begin{itemize}
    \item $h\A$ is a triangulated category;
    \item $\A_{o}$ has finite limits and finite colimits.
\end{itemize}
\end{remark}

Now, we tackle the idempotent complete condition and start by recalling a result about idempotent complete $\infty$-categories proved by Lurie.

\begin{proposition}[{\cite[Lemma 5.4.2.4]{HTT}}]
\label{propIdemComplCompact}
Let $\B$ be a small $\infty$-category. The following two conditions are equivalent:
\begin{itemize}
    \item the $\infty$-category $\B$ is idempotent complete, \Cref{exIdempo};
    \item the $\infty$-category $\B$ is equivalent to the compact $\infty$-subcategory of a $\omega$-accessible $\infty$-category.
\end{itemize}   
\end{proposition}

\begin{definition}
\label{defIndComplespectral}
We say that a spectral $\infty$-category $\A\in \Cat^{\Sp}$ is \textit{idempotent complete} if its underlying $\infty$-category $\A_o$ satisfies one of the equivalent conditions in \Cref{propIdemComplCompact}. 
\end{definition}

\begin{definition}
\label{defCategoryPerfect}
We say that a spectral $\infty$-category $\A$ is \textit{perfect} if it is small, idempotent complete, and pretriangulated. We denote by $\Cat^{\Sp,perf}$ the full $\infty$-subcategory of $\Cat^{\Sp}$ whose elements are small perfect spectral $\infty$-categories. We call $\Cat^{\Sp}$ the \textit{$\infty$-category of perfect spectral $\infty$-categories}. 
\end{definition}

\begin{remark}[Technical Point]
\label{remkTechinicalPoint2}
Let $\A$ be a perfect spectral $\infty$-category. It is important to emphasize that its underlying $\infty$-category is small.
\end{remark}

\begin{remark}
\label{rmkCatSpSpsonoCatEx}
It is interesting that $\Cat^{\Sp,perf}\subseteq \Cat^{\Sp}$ is a full $\infty$-subcategory and we do not request that the underlying $\infty$-functor $f_{o}:\A_o\to \B_{o}$ of a $\Sp$-enriched $\infty$-functor $f:\A\to \B$ is exact (i.e., it is the correct morphism between stable $\infty$-categories). This is because in \cite[Proposition 9.19]{heine2023equivalence}, the author proves that for any $\A$ and $\B$ pretriangulated spectral $\infty$-categories (in particular, that holds for perfect $\infty$-categories), the underlying $\infty$-functor $f_{o}:\A_{o}\to\B_{o}$ of an $\Sp$-enriched $\infty$-functor $f:\A\to\B$ is an exact $\infty$-functor. Moreover, he proves that the underlying $\infty$-category $\infty$-functor \[(-)_o:\Cat^{\Sp,perf}\to \Cat^{Ex}\] is conservative, \cite[Example 5.22]{DoniHigherCategorical}.
\end{remark}

\begin{definition}
\label{defcompactlygeneratedspectralcategory}
Let $\A$ be a pretriangulated spectral $\infty$-category. We say that it is compactly generated if it satisfies the following conditions:
\begin{itemize}
\item $\A_{o}$ admits small coproducts;
\item $\A_{o}$ is $\omega$-accessible;

\item for each object $X\in\A$, the $\infty$-functor $\A(-,X):\A^{op}_{o}\to \Sp$ preserves small products
\end{itemize}
\end{definition}

\begin{remark}[{\cite[Proposition 1.4.4.1]{HA}}]
\label{rmkPresentableEnrichedColimits}
Every stable $\infty$-category with small (co)products admits all small (co)limits, so $\A_o$ in \Cref{defcompactlygeneratedspectralcategory} is a compactly generated $\infty$-category. In particular, $\A_o$ is not usually small. For the same reason, $\A_{o}(-,X):\A^{op}\to \Sp$ preserves all small limits, 
and this implies that $\A$ has $\Sp$-enriched colimits.
Indeed, let $d:I\to \A_{o}$ be a diagram and let $colim_{i\in I}d\in \A_o$ be the nadir of the colimit (which is an apex of the limit of $d^{op}$ in $\A_o^{op}$, see \cite[Proposition 12.1.7]{RiehlElements}), so, for each $X\in\A_{o}$, the following equivalence holds:
\[\A(colim_{i\in I} d(i),X)\simeq \A(lim_{i\in I^{op}}d^{op}(i),X)\simeq lim_{i\in I^{op}}\A(d(i),X),\]
this is the $\infty$-categorical of the ordinary case, see \cite[\S A.5]{RiehlElements}. 
\end{remark}

Let $\A$ be a pretriangulated compactly generated spectral category. From the above remark, we know that $\A_o$ is a compactly generated stable $\infty$-category.
Just like being stable is the $\infty$-categorical analog of being triangulated for an ordinary category, being compactly generated for a stable $\infty$-category is the $\infty$-categorical analog of being a $\omega$-well generated triangulated category (or $\omega$-compactly generated triangulated category). If you are not familiar with what this concept is, you can find a presentation in \cite{canonaco2017tour} or in \cite{neeman2014triangulated}.

We recall the definition in \cite{neeman2014triangulated} of compact object in a triangulated category admitting small coproducts.

\begin{definition}
\label{defcompactNeeman}
Let $\mathcal{T}$ be a triangulated category admitting small coproducts. An object $c\in \mathcal{T}$ is said to be compact if the corepresentable functor $\Hom_{\C}(c,-):\C\to \Top$ commutes with small coproducts.
\end{definition}

\begin{proposition}
Let $\A$ be a stable $\infty$-category, then $\A$ is compactly generated if and only if:

\begin{itemize}
 
    \item $\A$ admits small coproducts;

    \item and the triangulated category $h\A$ is $\omega$-well generated.

\end{itemize}

In particular, an object of $\A$ is a compact object in the triangulated category with small coproduct $h\A$, \Cref{defcompactNeeman}, if and only if it is a compact object in the stable $\infty$-category admitting filtered colimits $\A_o$, \Cref{defcompactobjects}.    
\end{proposition}
\begin{proof}
This result is proven by Lurie in \cite[Proposition 1.4.4.1, Proposition 1.4.4.2]{HA}, \cite[Theorem 5.1.1.1]{HTT}.
\end{proof}

\begin{definition}
\label{defCompactlyGeneratedSpectraCategory}
We denote by $\Cat^{\Sp,h\text{-}proj}$ the $\infty$-subcategory of $\widehat{\Cat}^{\Sp}$ whose objects are compactly generated pretriangulated spectral $\infty$-categories and whose morphisms are $\Sp$-enriched $\infty$-functors such that their underlying $\infty$-functor preserves small coproducts. We call the objects of $\Cat^{\Sp,h\text{-}proj}$ the \textit{h-project spectral $\infty$-categories}. 
\end{definition}

\begin{definition}
\label{defSpectralCategoryofhprojectomega}
We denote by $\Cat^{\Sp,\omega}$ the $\infty$-subcategory of $\widehat{\Cat}^{\Sp}$ whose objects are compactly generated pretriangulated spectral $\infty$-categories and whose morphisms are $\Sp$-enriched $\infty$-functors such that their underlying $\infty$-functor preserves small coproducts and compact objects. We call the objects of $\Cat^{\Sp,\omega}$ the \textit{h-project spectral $\infty$-category}, because they are the same as $\Cat^{\Sp,h\text{-}proj}$. 
\end{definition}

\begin{definition}
\label{defpresentablespectralcategory}
Let $\A$ be a pretriangulated spectral $\infty$-category. We say that it is presentable if it satisfies the following conditions:
\begin{itemize}
\item $\A_o$ is a presentable $\infty$-category;
\item for each object $X\in\A$, the $\infty$-functor $\A(-,X):\A^{op}_{o}\to \Sp$ preserves all small limits.
\end{itemize}    
\end{definition}

\begin{definition}
\label{defcategorypresentablespectral}
We denote by $\Cat^{\Sp,L,St}$ the $\infty$-subcategory of $\widehat{\Cat}^{\Sp}$ whose objects are presentable pretriangulated spectral $\infty$-categories and whose morphisms are $\Sp$-enriched $\infty$-functors such that their underlying $\infty$-functor preserves small colimits. We call the objects of $\Cat^{\Sp,L,St}$ the \textit{presentable pretriangulated spectral $\infty$-categories}. 
\end{definition}

The pretriangulated presentable spectral $\infty$-categories have an important property.

\begin{proposition}
\label{propPretriangulatedCompactlyGeneratedareTensore}
Let $\A$ be a pretriangulated presentable spectral $\infty$-category, then $\A$ is a tensored $\Sp$-enriched $\infty$-category.
Vice versa, every tensored $\Sp$-enriched $\infty$-category $\A$ whose underlying $\infty$-category $\A_o$ is presentable is a pretriangulated presentable spectral $\infty$-category.     
\end{proposition}
\begin{proof}
 By definition, the underlying $\infty$-category $\A_o$ is presentable and stable so it is a left $\Sp$-module object in $Pr^{L}$. In particular, there is an $\infty$-functor \[-\otimes- :\Sp\times \A_o\to \A_o\] which preserves small colimits componentwise. It is straightforward to prove that this $\infty$-functor satisfies the properties in \cite[Definition 5.32]{DoniHigherCategorical}. The vice versa is obvious.
\end{proof}

\begin{remark}
\label{rmkCatSpLStPrSp}
As a consequence of \Cref{propPretriangulatedCompactlyGeneratedareTensore}, the $\infty$-category $\Cat^{\Sp, L,St}$ is the $\infty$-category $Pr^L_{\Sp}$ defined in \cite{heine2023equivalence}. Indeed, $Pr^L_{\Sp}$ is the $\infty$-subcategory of $\widehat{\Cat}^{\Sp}$ generated by tensored $\Sp$-enriched $\infty$-categories whose underlying $\infty$-category is presentable and whose morphisms are the same as $\Cat^{\Sp,L,St}$. In particular, in \cite[Theorem 1.2]{heine2023equivalence}, the author proves the fundamental result for us that \[\Cat^{\Sp, L,St}=Pr^L_{\Sp}\simeq \mathrm{LMod}_{\Sp}(Pr^{L})\simeq Pr^{L,St}.\]
\end{remark}

Now we add a $k$-linearization to the above definitions. The idea of the $k$-linearization definitions is very simple: the properties inside the underlying $\infty$-category $\A_{o}$ remain unchanged and the properties associated with the (co)limits-preservation $\infty$-functors remain the same but the $\infty$-functors will target the $\infty$-category $\LModk$ instead of $\Sp$. We will find out that even this difference can be eliminated.

\begin{definition}
\label{defPretriangulateddgcategories}
Let $\A$ be a dg-category, we will say that it is \textit{pretriangulated} if its underlying $\infty$-category $\A_{o}$ is stable
and, for each objects $Y\in\A_{o}$, the $\infty$-functors $\A(Y,-):\A_{o}\to \LModk$ and $\A(-,Y):\A_{o}^{op}\to \LModk$ are exact.
\end{definition}

\begin{definition}
We denote by $\Cat^{\LModk,St}$ the full $\infty$-subcategory of $\Cat^{\LModk}$ whose objects are small pretriangulated dg-categories.
\end{definition}

\begin{remark}
The same deduction of \Cref{rmkEnrichedcolimitsEnrichedLimits} holds for pretriangulated dg-category, with the only difference that it admits finite $\LModk$-enriched colimits and finite $\LModk$-enriched limits, instead of $\Sp$-enriched.
\end{remark}

\begin{definition}
\label{defIndCompledgcat}
We say that a dg-category $\A\in \Cat^{\LModk}$ is \textit{idempotent complete} if its underlying $\infty$-category $\A_o$ satisfies one of the equivalent conditions in \Cref{propIdemComplCompact}. 
\end{definition}

\begin{definition}
\label{defdgCategoryPerfect}
We say that a dg-category $\A$ is \textit{perfect} if it is small, idempotent complete, and pretriangulated. Moreover, we denote by $\Cat^{\LModk,perf}$ the full $\infty$-subcategory of $\Cat^{\LModk}$ whose elements are small perfect dg-categories. We call $\Cat^{\LModk,perf}$ the \textit{$\infty$-category of small perfect dg-categories $\infty$-categories}. 
\end{definition}


\begin{definition}
\label{defcompactlygenerateddgcategory}
Let $\A$ be a pretriangulated dg-category, we say that it is compactly generated if it satisfies the following conditions:
\begin{itemize}
\item $\A_{o}$ admits small coproducts,
\item  $\A_{o}$ is $\omega$-accessible;
\item for each object $X\in\A$, the $\infty$-functor $\A(-,X):\A^{op}\to \LModk$ preserves small products.
\end{itemize}
\end{definition}

\begin{definition}
\label{defCompactlygenerateddg-category}
We denote by $\Cat^{\LModk,h\text{-}proj}$ the $\infty$-subcategory of $\widehat{\Cat}^{\LModk}$ whose objects are compactly generated pretriangulated dg-categories and whose morphisms are $\LModk$-enriched $\infty$-functors such that their underlying $\infty$-functor preserves small coproducts.  
We call the objects of $\Cat^{\LModk,\omega}$ \textit{h-project dg-categories}. 
\end{definition}

\begin{remark}
Let $\A$ be a pretriangulated compactly generated dg-category.
Since $\A_o$ is stable and has small coproducts, then it has small colimits, \cite[Proposition 1.4.4.1]{HA}. In particular, $\A_o$ is a compactly generated $\infty$-category.
For the same reason, the morphisms in $\Cat^{\LModk,h\text{-}proj}$ preserve all small colimits. In fact, they preserve small coproducts by definition and they are exact, see \Cref{rmkCatSpSpsonoCatEx}. 
\end{remark}

\begin{definition}
\label{defdgCategoryofhomega}
We denote by $\Cat^{\LModk,\omega}$ the $\infty$-subcategory of $\widehat{\Cat}^{\LModk}$ whose objects are compactly generated pretriangulated dg-categories and whose morphisms are $\LModk$-enriched $\infty$-functors such that their underlying $\infty$-functor preserves small coproducts and compact objects. 
We call the objects of $\Cat^{\LModk,\omega}$ \textit{h-project dg-categories}. 
\end{definition}

\begin{definition}
\label{defpresentabledgcategory}
Let $\A$ be a pretriangulated dg-category,
we say that it is presentable if it satisfies the following conditions:
\begin{itemize}
\item $\A_o$ is a presentable $\infty$-category;
\item for each object $X\in\A$, the $\infty$-functor $\A(-,X):\A^{op}\to \LModk$ preserves all small limits.
\end{itemize}    
\end{definition}

\begin{remark}
The same deduction of \Cref{rmkPresentableEnrichedColimits} holds for a compactly generated pretriangulated dg-category and for a presentable pretriangulated dg-category, with the only difference that it admits small $\LModk$-enriched colimits, instead of $\Sp$-enriched.
\end{remark}

\begin{definition}
\label{defdgCategoryofhproject}
We denote by $\Cat^{\LModk,L,St}$ the $\infty$-subcategory of $\widehat{\Cat}^{\LModk}$ whose objects are presentable pretriangulated dg-categories and whose morphisms are $\LModk$-enriched $\infty$-functors such that their underlying $\infty$-functor preserves small colimits. 
We call the objects of $\Cat^{\Sp,L,St}$ \textit{presentable pretriangulated dg-categories}. 
\end{definition}

The pretriangulated presentable dg-categories have an important property.

\begin{proposition}
\label{propPretriangulatedCompactlyGeneratedareTensoreLMODK}
Let $\A$ be a pretriangulated presentable dg-category, then $\A$ is a tensored $\LModk$-enriched $\infty$-category.
Vice versa, every tensored $\LModk$-enriched $\infty$-category $\A$ whose underlying $\infty$-category $\A_o$ is presentable is a pretriangulated presentable spectral $\infty$-category.     
\end{proposition}
\begin{proof}
The proof is the same as \Cref{propPretriangulatedCompactlyGeneratedareTensore}.
\end{proof}

\begin{remark}
\label{rmksuperimportant}
As a consequence of \Cref{propPretriangulatedCompactlyGeneratedareTensoreLMODK}, $\Cat^{\LModk, L,St}$ is the $\infty$-category $Pr^L_{\LModk}$ defined in \cite{heine2023equivalence}. Indeed, $Pr^L_{\LModk}$ is the $\infty$-subcategory of $\widehat{\Cat}^{\LModk}$ generated by tensored $\LModk$-enriched $\infty$-categories whose underlying $\infty$-category is presentable and whose morphisms are the same as $\Cat^{\LModk,L,St}$. In particular, in \cite[Theorem 1.2]{heine2023equivalence}, the author proves the fundamental result for us that \begin{equation}
    \Cat^{\LModk, L,St}= Pr^L_{\LModk}\simeq \mathrm{LMod}_{\LModk}(Pr^{L}).
    \end{equation}
\end{remark}

Let us recall some results that allow us to reduce the study of dg-categories to the study of their underlying spectral $\infty$-categories $\A_{sp}$.

\begin{remark}[{\cite[Corollary 4.2.3.5.]{HA}\cite[Corollary 4.2.3.3]{HA}}]  
\label{rmkforgetfulcreatescolimitsandlimits}
Let $U_k:\LModk\to \Sp$ the forgetful $\infty$-functor: 
\begin{itemize}
\item $U_k$ preserves, reflects, and creates small colimits;  
\item $U_k$ preserves and reflects small limits.
\end{itemize}

In \cite[Corollary 6.9]{DoniHigherCategorical}, we proved that there exist an equivalent to $\Cat^{\LModk}\simeq \mathrm{LMod}_{\Hk}(\Cat^{\Sp})$.
As above, the forgetful $\infty$-functor \[\mathcal{U}_k:\Cat^{\LModk}\simeq \mathrm{LMod}_{\Hk}(\Cat^{\Sp})\to \Cat^{\Sp}\] has the properties: 
\begin{itemize}
\item $\mathcal{U}_k$ preserves, reflects, and creates small colimits;  
\item $\mathcal{U}_k$ preserves and reflects small limits.
\end{itemize}

Note that the forgetful $\infty$-functor $\mathcal{U}_k$ can be described in at least two other ways: 
as the evaluation $\infty$-functor \[ev_{*}:\Cat^{\LModk}\simeq Fun^{Alg(\Sp)}(\un{\Hk},\Cat^{\Sp})\to \Cat^{\Sp},\] or as the change of base $\infty$-functor $U_{!}$, \cite[Exampe 5.25]{DoniHigherCategorical}. 
\end{remark}

Using \Cref{rmkforgetfulcreatescolimitsandlimits}, we can reformulate the above definitions.

\begin{proposition}
 \label{defPretriangulateddgcategories2}
Let $\A$ be a dg-category. Then the following properties are equivalent:
\begin{itemize}
    \item $\A$ is pretringulated;  
    \item its underlying $\infty$-category $\A_{o}$ is stable and for each object $Y\in\A_{o}$, the $\infty$-functors $\mathcal{U}_k\circ \A(Y,-):\A_{o}\to \LModk\to \Sp$ and $\mathcal{U}_k\circ\A(-,Y):\A_{o}^{op}\to \LModk\to \Sp$ are exact;
   \item a dg-category $\A\in \Cat^{\LModk}$ is pretriangulated if its underlying spectral $\infty$-category $\A_{sp}$ is. 
\end{itemize}
\end{proposition}
\begin{proof}
    Follows from \Cref{rmkforgetfulcreatescolimitsandlimits}.
\end{proof}

\begin{proposition}
\label{propNewDefinitionPretrdgcatIdemp}
 Let  $\A$ be a pretrinagulated $\infty$-category, then the following properties are equivalent:
 \begin{itemize}
     \item it is idempotent complete,
     \item its underlying spectral $\infty$-category $\A_{sp}$ is idempotent complete;
     \item its underlying $\infty$-category $\A_{o}$ is idempotent complete.
 \end{itemize}
\end{proposition}

\begin{proof}
    Follows from \Cref{rmkforgetfulcreatescolimitsandlimits}.
\end{proof}

\begin{remark}
\label{CatSpperfisapullback}
From \Cref{defPretriangulateddgcategories2} and \Cref{propNewDefinitionPretrdgcatIdemp}, we obtain that $\Cat^{\LModk,perf}$ is the pullback square in the following diagram:

\begin{equation}
\label{eqCatPerfModkisPullback}
\begin{tikzcd}[ampersand replacement=\&]
	{\Cat^{\LModk,perf}} \& {\Cat^{\LModk}} \\
	{\Cat^{\Sp,perf}} \& {\Cat^{\Sp}}
	\arrow["\mathcal{U}_{k}^{perf}"',from=1-1, to=2-1]
	\arrow[from=2-1, to=2-2]
	\arrow[from=1-1, to=1-2]
	\arrow["\mathcal{U}_{k}",from=1-2, to=2-2]
	\arrow["\lrcorner"{anchor=center, pos=0.125}, draw=none, from=1-1, to=2-2]
\end{tikzcd}
\end{equation}
Since the class of conservative $\infty$-functors is the right class of a factorization system in $\Cat$, see \cite[Example 3.1.7 (f)]{anel2022left}, the class of conservative $\infty$-functors is stable under pullback (or base change). This implies that the $\infty$-functor $\mathcal{U}_{k}^{perf}$ in \eqref{eqCatPerfModkisPullback} is conservative.
So we deduce that the $\infty$-functor \[\Cat^{\LModk,perf}\xrightarrow{\mathcal{U}_k^{perf}}\Cat^{\Sp,perf}\xrightarrow{(-)_o} \Cat^{Ex}\] is conservative because it is a composition of conservative $\infty$-functors, see \Cref{rmkCatSpSpsonoCatEx}.
\end{remark}

\begin{proposition}
\label{propnewsdefinitionhprojdgcat}
Let $\A$ be a pretriangulated dg-category, then the following properties are equivalent:
\begin{itemize} 
    \item $\A$ is compactly generated;

    \item its underlying $\infty$-category $\A_{o}$ is compactly generated and for each object $Y\in\A_{o}$, the $\infty$-functor $\mathcal{U}_k\circ\A(-,Y):\A_{o}^{op}\to \LModk\to \Sp$ preserves small limits;

    \item the underlying spectral $\infty$-category $\A_{sp}$ is compactly generated;
\end{itemize}    
\end{proposition}
\begin{proof}
    Follows from \Cref{rmkforgetfulcreatescolimitsandlimits}.
\end{proof}

\begin{remark}
\label{CatSpgprojisapullback}
From \Cref{propnewsdefinitionhprojdgcat}, we obtain that $\Cat^{\LModk,\omega}$ is the pullback square in the following diagram:

\[\begin{tikzcd}[ampersand replacement=\&]
	{\Cat^{\LModk,\omega}} \& {\Cat^{\LModk}} \\
	{\Cat^{\Sp,\omega}} \& {\Cat^{\Sp}}
	\arrow[from=1-1, to=2-1]
	\arrow[from=2-1, to=2-2]
	\arrow[from=1-1, to=1-2]
	\arrow[from=1-2, to=2-2]
	\arrow["\lrcorner"{anchor=center, pos=0.125}, draw=none, from=1-1, to=2-2]
\end{tikzcd}\]
\end{remark}

\begin{proposition}
    Let $\A$ be a pretriangulated compactly generated dg-category, then $\A$ is a tensored $\LModk$-enriched $\infty$-category.    
\end{proposition}
\begin{proof}
 This result is a consequence of \Cref{propPretriangulatedCompactlyGeneratedareTensore}.
\end{proof}

Now we prove the $\Sp$-enriched version and the $\LModk$-enriched version of (\ref{eqCatPerfPrStLOmega}).


\begin{proposition}
\label{propCatSpPerfPrOmega}
There exists an equivalence of $\infty$-categories 
\begin{equation}
    \label{eqimportant}
    M: \Cat^{\Sp,perf} \simeq \Cat^{\Sp,\omega}: Y.
\end{equation}
\end{proposition}

\begin{proof}
In \cite[Theorem 1.2]{heine2023equivalence}, the author proves that there exists an equivalence $\Cat^{\Sp,St,L} \simeq Pr^{L,St}$. We restrict this equivalence via the inclusion $Pr^{L,St}_{\omega} \to Pr^{L,St}$, and we obtain an equivalence:
\begin{equation}
\label{eq1edo}
\Cat^{\Sp,\omega} \simeq Pr^{L,St}_{\omega}.
\end{equation}
In the same paper \cite[Theorem 9.18]{heine2023equivalence}, the author proves that there exists an equivalence $\Cat^{\Sp,St} \simeq \Cat^{Ex}$. We restrict this equivalence via the inclusion $\Cat^{perf} \to \Cat^{Ex}$, and we obtain an equivalence:  
\begin{equation}
\label{eqedo2}
\Cat^{\Sp,perf} \simeq \Cat^{perf}.
\end{equation}

So we obtain the thesis from (\ref{eqedo2}) and (\ref{eq1edo}) and (\ref{eqCatPerfPrStLOmega}).
\end{proof}

\begin{remark}
\label{rmkCatSpPerfSmall}
The above proposition is important because the elements in $\Cat^{\Sp,\omega}$ are by definition not necessarily small; instead, the elements of $\Cat^{\Sp,perf}$ are small.
\end{remark}

\begin{remark}
\label{rmkDescriptionofequivalenceCatspperfCatperf}
The statement of \ref{propCatSpPerfPrOmega} is interesting. 

Combining it with a few results from \cite{heine2023equivalence} and in \cite{blumberg2013universal}, we can deduce some important facts. Indeed, the $\infty$-functor $M:\Cat^{\Sp,perf}\to Cat^{\Sp,L,St}_{\infty,\omega}$ in (\ref{eqimportant}) postcomposed with the $\infty$-functor underlying $\infty$-category \[M(-)_o := \Cat^{\Sp,perf} \to Cat^{\Sp,\omega}_{\infty} \xrightarrow{(-)_{o}} Pr_{\omega}^{L},\] can be described on the objects in three ways.

Let $\A$ be an object of $\Cat^{\Sp,perf}$:
\begin{itemize}
    \item[(i)] $M(\A)_{o}$ is the ind-completion of the underlying category of $\A$, i.e., $M(\A)_{o} \simeq Ind(\A_o)$;
    \item[(ii)] $M(\A)_{o}$ is the $\infty$-category of $\Sp$-enriched $\infty$-functors $P_{\Sp}(\A)$. In particular, $M(\A)_{o}$ is the Hinich's Yoneda Embedding;
    \item[(iii)] $M(\A)_o$ is the $\infty$-category of exact $\infty$-functors $Fun^{Ex}(\A_o,\Sp)$.
\end{itemize}
Heine proves that $(i)$ is equivalent to $(ii)$, whereas Tabuada-Blumberg-Gepner prove that $(i)$ is equivalent to $(iii)$.

Moreover, if we further compose $M(-)_o$ with the equivalence $(-)_{\omega}:Pr^{L,St}_{\omega}\to \Cat^{perf}$, we obtain an $\infty$-functor 
\[M(-)_{o,\omega}:\Cat^{\Sp,perf}\to Pr_{\omega}^{L} \xrightarrow{(-)_{\omega}} \Cat^{perf}.\]
Also, $M(-)_{o,\omega}$ can be described on the objects in three ways. Let $\A$ be an object of $\Cat^{\Sp,perf}$:
\begin{itemize}
    \item[(1)] $M(\A)_{o,\omega}$ is the full $\infty$-subcategory generated by the compact objects of $Ind(\A_{o})$, i.e., $M(\A)_{o,\omega} \simeq Ind(\A_o)^{\omega}$. But by definition, $\A_o \simeq Ind(\A_o)^{\omega}$, so $M(\A)_{o,\omega}$ is the underlying $\infty$-category;
    \item[(2)] $M(\A)_{o,\omega}$ is the full $\infty$-subcategory generated by the compact objects of $\Sp$-enriched $\infty$-functors $\infty$-category $P_{\Sp}(\A)$, i.e., $M(\A)_{o,\omega} \simeq P_{\Sp}(\A)^{\omega}$;
    \item[(3)] $M(\A)_{o,\omega}$ is the full $\infty$-subcategory generated by the compact objects of $Fun^{Ex}(\A_o,\Sp)$, i.e., $M(\A)_{o,\omega} \simeq Fun^{Ex}(\A_o,\Sp)^{\omega}$.
\end{itemize}

Instead, if we post-compose the $\infty$-functor $Y:\Cat^{\Sp,\omega} \simeq \Cat^{\Sp,perf}$ with the underlying $\infty$-category functor, we obtain an $\infty$-functor:
\[Y(-)_o:\Cat^{\Sp,\omega} \simeq \Cat^{\Sp,perf} \xrightarrow{(-)_o} \Cat^{perf},\]
which carries an object $\A$ of $\Cat^{\Sp,\omega}$ in the full $\infty$-subcategory generated by the compact objects of its underlying $\infty$-category, i.e., $Y(\A)_o \simeq (\A_o)^{\omega}$. 

We conclude that the adjoint equivalence $M \dashv Y$ is the $\Sp$-enriched analog of \eqref{eqCatPerfPrStLOmega}. 
\end{remark}

Now, we are going to $k$-linearize \Cref{propCatSpPerfPrOmega}.

\begin{corollary}
\label{CorCatLModPerfCatLModhprojequivalence}
There exists an equivalence of $\infty$-categories
\[ \Cat^{\LModk,perf} \simeq \Cat^{\LModk,\omega}. \]
\end{corollary}

\begin{proof}
Since there is the equivalence between $\Cat^{\Sp,perf} \simeq \Cat^{\Sp,\omega}$, by \Cref{propCatSpPerfPrOmega}, they are apex of the same pullback diagram, see \Cref{CatSpperfisapullback} and \Cref{CatSpgprojisapullback}.
So, by the universal property of the pullback, we obtain the thesis.
\end{proof}

\begin{remark}
\label{rmkSize}
The result in \Cref{CorCatLModPerfCatLModhprojequivalence} is important because it allows us to work with presentable $\infty$-categories using small $\infty$-categories; i.e., the objects in $\Cat^{\LModk,\omega}$ are compactly generated and not necessarily small, whereas the objects in $\Cat^{\LModk,perf}$ are small.
\end{remark}

\section{$k$-linear stable $\infty$-categories}
\label{secklinearstableCat}

\begin{warning}
\label{WarningDEfKlinear}
In \cite{SAG}, the author chooses that a $k$-linear (stable) $\infty$-category is always a presentable $\infty$-category. While for us, a $k$-linear (stable) $\infty$-category is not necessarily presentable. Here, we do not define what a $k$-linear stable $\infty$-category is, and we hope that this topic is for future work.
\end{warning}

In this section, we will define what a $k$-linearization is for a general $\infty$-category (not necessarily enriched) and analyze some specific cases.

The idea is to use the construction and results seen in \cite{DoniHigherCategorical} and to remove from them the dependence on enrichment. In the work just mentioned, we have proved that an $\Hk$-linearization of an $\infty$-category $\A$, i.e., $\A$ with an $\LModk$-enrichment, is equivalent to giving a left $\Hk$-module of $\Cat^{\Sp}$, i.e., a pair $(\hat{A}, \un{\Hk}\otimes^{\Sp}\hat{A}\to\hat{A})$ where the first component is a spectral $\infty$-category (or a $\Sp$-enriched $\infty$-category) and the second component is a left $\Hk$-action on $\hat{A}$.

In turn, in \cite[Corollary 6.10 and Corollary 6.12]{DoniHigherCategorical}, we have proved that a left $\Hk$-module can be seen as an enriched $\infty$-functor from $\un{\Hk}$ to an Internal-hom, $\un{\Hk}\to \un{(\hat{A})}^{\un{(\hat{A})}}$.

We deduce that a general $k$-linearization (probably too general to be useful) is a left $\un{\Hk}$-module in $\Cat$ that can be seen as $\infty$-functor from $\un{\Hk}$ to the internal-hom of $\Cat$ $A:\un{\Hk}\to Fun(-,-)$, see \Cref{defKlinear}. 

Some $\Sp$-enrichement can be internalized in the theory of $\infty$-categories, in particular, it is uniquely determined by limits and colimits in $\A_o$. Something very similar happens in the ordinary category theory too: additive $\infty$-categories and abelian categories have the enrichment which is uniquely determined by their limits and their colimits. We will see this fact in \Cref{subsectionPresentableStableKlinearCat}.
So, the $\infty$-functor $A$ must have the property that \enquote{If $\A_{o}$ has suitable colimits or/and limits and $A$ respects these limits and colimits, then it defines uniquely a $k$-linearization.}

There is a drawback in our last sentence: what is an $\infty$-functor $\un{\Hk}\to Fun(-,-)$? Or equivalently, what is a left $\un{\Hk}$-module in $\Cat$? an $\infty$-functor $\un{\Hk}\to Fun(-,-)$ is not well defined because the source and the target live in two different $\infty$-categories: the former lives in $\Cat^{\Sp}$ while the latter lives in $\Cat$.

In \Cref{subsecGeneraleklinear}, we will solve the drawback using a different point of view of spectra and explain how to formalize \enquote{general $k$-linearization}. Probably this generalization is useless, but we think that it can help to understand the reasoning behind this idea. Instead, some less general $k$-linearizations can be interesting and we hope that they will be the topic of some future work.

In \Cref{subsectionPresentableStableKlinearCat}, we will analyze the presentable and stable case: here we will find that our definition and Lurie's definition are equivalent. This case has the drawback that the presentable $\infty$-categories are big and might not be the correct notion for our main theorem. That's why there is the \Cref{subsecMoriaTheoryklinear} where we will define the correct setting for our Main Theorem.





\subsection{Defining $k$-linearization for general $\infty$-categories}
\label{subsecGeneraleklinear}

This subsection provides informative content, and we hope to return in the future to analyze the ideas it contains; thus, not all the ideas are proven.

First, we recall that a connective spectrum (e.g., $\Hk$) can be seen as a grouplike commutative monoid in $\mathcal{S}$, \cite[Remark 5.2.6.26 and \S 1.4 (D)]{HA}. In the following remark, we summarize this viewpoint.

\begin{remark}[{\cite[\S 0.2: Spectra are generalized abelian groups]{SAG}}]
\label{RmkSpectraAreGeneralizedAbelianGroups}
Let $E$ be a spectrum. Then, the $0^{th}$ space $\Omega^{\infty}E$ is an example of an $\E$-space: that is, it can be equipped with an addition law
\[+:\Omega^{\infty}E\times\Omega^{\infty}E\to \Omega^{\infty}E\]
that is unital, commutative, and associative up to coherent homotopy. Moreover, the construction $E\to \Omega^{\infty}E$ restricts to an equivalence from the $\infty$-category $Sp^{cn}$ of connective spectra to the $\infty$-category $CAlg^{gp}(\Top)$ of grouplike $\E$-spaces
\[\Omega^{\infty}:\Sp^{cn}\simeq CAlg^{gr}(\Top).\]
An $\E$-space $A$ is said to be grouplike if the addition on $A$ exhibits the set of connected components $\pi_0A$ as an abelian group.
\end{remark}

This is the spectral algebraic view of a spectrum, and it is the fundamental observation that allows us to define a $k$-linearization without a starting $\Z$-linearization.

In particular, we can see $\Hk$ as an $\infty$-category as a space via the forgetful $\infty$-functor $U:CAlg^{gp}(\Top)\to \Top$, which is monoidal with respect to the cartesian structure of $\Top$. Thus, it induces a morphism between algebras:
\begin{equation}
U:Alg(CAlg^{gp}(\Top))\to Alg(\Top);
\end{equation}
and between commutative algebras:
\begin{equation}
\label{eqForgetfulMongr}
Alg(U):Alg(CAlg^{gp}(\Top))\to Alg(\Top).
\end{equation}

$Alg(U)$ allows us to see an $\Hk$ only as commutative algebras in the $\infty$-category $\Top$. Postcomposing the $\infty$-functor in \eqref{eqForgetfulMongr} with the underlying $\infty$-functor $(-)$ (it is defined in \cite{DoniHigherCategorical}), we obtain the next definition. 

\begin{definition}
\label{defMannaggia}
We call the \textit{overline $\infty$-functor} the following $\infty$-functor:
    \[\ov{(-)}:= \un{(-)}\circ Alg(U):Alg(CAlg^{gp}(\Top))\to Alg(\Top)\subseteq \Cat^{\Cat}:\]
    \[A\mapsto \ov{\A},\]
and we denote it by $\ov{(-)}$.
\end{definition}

It is well-known that $(\Cat,\times,*, Fun(-,-))$ is a closed monoidal $\infty$-category. This implies that $\Cat$ is a $\Cat$-enriched $\infty$-category.
Since $\Cat$ is a presentably monoidal $\infty$-category, the $\infty$-category of $\Cat$-enriched $\infty$-categories $\Cat^{\Cat}$ has a monoidal structure whose tensor is the $\Cat$-enriched tensor $-\otimes^{\Cat}-$ and whose unit is $\un{\mathbb{1}}$.
Furthermore, $\Cat$ has a canonical structure of left $\Cat$-module of $Pr^{L}$ and, since the full inclusion $\Top\to \Cat $ is monoidal, it is also a left $\Top$-module of $Pr^{L}$.

In \cite[Corollary 6.10]{DoniHigherCategorical}, we have seen that a $k$-linearization is a $\Cat^{\Sp}$-enriched $\infty$-functor with source $\2un{\Hk}$ and target $\Cat^{\Sp}$-enriched $\infty$-categories $\Cat^{\Sp}$: note that the objects of $\Cat^{\Sp}$ are non-geometrical $\infty$-categories.
In short, such an $\infty$-functor is a way to linearize a non-geometrical object of $\Cat^{\Sp}$. 
Now, we have the idea and the tools to give the first definition of a general $k$-linearization. 

\begin{definition}
\label{defKlinear}
We denote by $LinCat_{\Hk}$ the $\infty$-category 
$Fun^{\Cat}(\ov{\Hk},\Cat)$ and we call it the $\infty$-category of \textit{general $k$-linear $\infty$-categories}.
Moreover, we call $k$-linear $\infty$-categories the objects of $Fun^{\Cat}(\ov{\Hk},\Cat)$.
\end{definition}

\begin{remark}
\label{rmkInutile}
Since $\ov{\Hk}$ is an $\infty$-category, $\Cat$-enriched $\infty$-functors with source $\ov{\Hk}$ are actually $\infty$-functors. So we immediately have a redefinition of the $\infty$-category of general $k$-linear $\infty$-categories:
\[LinCat_{\Hk}=Fun^{\Cat}(\ov{\Hk},\Cat)\simeq Fun(\ov{\Hk},\Cat).\]
\end{remark}

Just as in \cite[Corollary 6.10]{DoniHigherCategorical}, we can immediately give another definition. 
\begin{proposition}
  There is an equivalence
  \[LinCat_{\Hk}\simeq \mathrm{LMod}_{\Hk}(\Cat).\]
\end{proposition}
\begin{proof}
Since the space of objects of $\ov{\Hk}$ is trivial, i.e., $\ov{\Hk}_{ob}\simeq *$, the equivalence follows by the definition of $\infty$-category of enriched $\infty$-functors, \cite{DoniHigherCategorical}.
\end{proof}

These definitions above are useless because they have lost the essential fact that a connective spectrum, as $\Hk$, is a generalized abelian group, \Cref{RmkSpectraAreGeneralizedAbelianGroups}.
We have lost this information when we use the forgetful $\infty$-functor (\ref{eqForgetfulMongr}) in the definition of $\ov{(-)}$. 

Instead of $\ov{(-)}$, we have to use the underline $\infty$-functor \[\un{(-)}:Alg(CAlg^{gr}(\Top))\to \Cat^{CAlg^{gr}(\Top)}.\] 

In particular, the $\infty$-category of $k$-linear $\infty$-categories should be an $\infty$-category of $CAlg^{gr}(\Top)$-enriched $\infty$-functors. 

In short, we have to restrict our attention to a suitable (not full) $\infty$-subcategory $S$ of $\Cat$ with the following features:
an $\infty$-category $\C$ belongs to $S$ if:
\begin{itemize}
    \item for each pair of objects $x,y\in S$, $\Hom_{S}(x,y)$ belongs to $CAlg^{gr}(\Top)\simeq \Sp^{cn}$, i.e., it is, at least, a $Sp^{cn}$-enriched $\infty$-category; 
    
    \item $\C$ admits $Y$-colimits, for a suitable collection of simplicial sets $Y$, such that the colimits \enquote{are well intertwined} with the structure of $CAlg^{gr}(\Top)$ of the Hom-spaces.  
\end{itemize}
An $\infty$-functor $F:\A\to \B$ belongs to $S$ if:  
\begin{itemize}
    \item for each pair of objects $x,y\in S$, $F:\Hom_{S}(x,y)\to \Hom(Fx,FY)$ belongs to $CAlg(\Top)$; 
    
    \item $F$ preserves $Y$-colimits.
\end{itemize}

Since $CAlg^{gr}(\Top)\simeq \Sp^{cn}$ has a monoidal structure \cite[Lemma 7.1.1.7.]{HA}, the definition of $\infty$-category of $k$-linear $S$ $\infty$-categories is 

\begin{definition}
\label{eqdefinitionideaDefKLinea1}
We denote by $LinCat^{S}_{\Hk}$ the $\infty$-category 
$Fun^{CAlg^{gr}(\Top)}(\un{Hk},S)\simeq Fun^{\Sp^{cn}}(\un{\Hk},S)$ and we call it the $\infty$-category of \textit{$S$ $k$-linear $\infty$-categories}.
Moreover, we call $S$ $k$-linear $\infty$-categories the objects of $LinCat^{S}_{\Hk}$.    
\end{definition}

The above definition is well posed if the $\infty$-category $Fun^{CAlg^{gr}(\Top)}(\Hk,S)$ is well defined; e.g., if $S$ is a left $Sp^{cn}$-module in $\Cat(\K)$.

In \cite[Lemma 7.1.1.7.]{HA}, the author proves that the inclusion $i:\Sp^{cn}\to Alg(\Sp)$ is monoidal. Tautologically, giving a $Sp^{cn}$-enriched $\infty$-functor with source $\un{\Hk}$ is the same as giving a $Alg(\Sp)$-enriched $\infty$-functor with source $\un{\Hk}$, so the next definition is the same as \Cref{eqdefinitionideaDefKLinea1} . 

\begin{definition}
\label{eqdefinitionideaDefKLinear}
We denote by $LinCat^{S}_{\Hk}$ the $\infty$-category 
$Fun^{Alg(\Sp)}(\un{\Hk},S)$ and we call it the $\infty$-category of \textit{$S$ $k$-linear $\infty$-categories}.
Moreover, we call $S$ $k$-linear $\infty$-categories the objects of $LinCat^{S}_{\Hk}$.    
\end{definition}
The above definition is well posed if the $\infty$-category $Fun^{\Sp}(\un{\Hk},S)$ is well defined; e.g., if $S$ is a left $Alg(Sp)$-module in $\Cat(\K)$.

In \Cref{subsectionPresentableStableKlinearCat} and \Cref{subsecMoriaTheoryklinear}, we will study the situation when $S$ is $ Pr^{L,St}$, $Pr^{L,St}_{\omega}$ or $\Cat^{perf}$. 
There are other possible choices for $S$: e.g. the $\infty$-category of additive $\infty$-categories $\Cat^{Add}$, the $\infty$-category of Grothendieck prestable $\infty$-categories $Groth$.
In this paper, we are not interested in these other choices, but we hope that they will be the subject of future work.



\subsection{ Presentable and stable $k$-linearity $\infty$-categories }
\label{subsectionPresentableStableKlinearCat}

In this subsection, we will define what a $k$-linearization is for a presentable stable $\infty$-category.

\begin{remark}
\label{rmkHKAlgebraDiAlgebra}
Since $\Hk$ is an object of $CAlg(\Sp)$, it is particularly an algebra in the symmetric monoidal $\infty$-category $Alg(\Sp)$. There exists the canonical forgetful $\infty$-functor:
\begin{equation}
\label{eqHKAlgebraDiAlgebra}
U:CAlg(\Sp)\simeq Alg_{\E}(\Sp)\to Alg_{\mathbb{E}_{2}}(\Sp)\simeq Alg(Alg(\Sp)).
\end{equation}
Also, in this setting, we have the underline $\infty$-functor:
\begin{equation}
    \un{(-)}:Alg(Alg(\Sp))\to \Cat^{Alg(\Sp)}.
\end{equation}
As a consequence, \Cref{propActionklinearcategory} and \Cref{propNewDefinition} will be well-posed.
This observation emphasizes why it is important that $\Hk$ is a commutative algebra; i.e., an $\E$-ring.
Actually, it is sufficient that it is an $\mathbb{E}_2$-ring, but in this article, we are not interested in this generalization.
\end{remark}

We know that the $\infty$-category $Pr^{L,St}$ is a left $Alg(\Sp)$-module object of $Pr^{L}$, see \eqref{eqTensorPrL}, so using the idea in \Cref{eqdefinitionideaDefKLinear}, the next definition is the natural one.

\begin{definition}
\label{defKlinearstablePresentable}
We denote by $LinCat_{\Hk}^{L,St}$ the $\infty$-category $Fun^{Alg(\Sp)}(\un{\Hk}, Pr^{L,St})$. Moreover, we call \textit{$k$-linear stable presentable $\infty$-categories} the objects of $LinCat_{\Hk}^{L,St}$, and we call $LinCat_{\Hk}^{L,St}$ the $\infty$-category \textit{of stable presentable k-linear} $\infty$-categories. 
\end{definition}

\begin{proposition}
\label{propActionklinearcategory}
There exists an equivalence of $\infty$-categories:
\[ LinFun^{L,St}_{\Hk}= Fun^{Alg(\Sp)}(\Hk,Pr^{L,St})\simeq \mathrm{LMod}_{\Hk}(Pr^{L,St}).\]
\end{proposition}
\begin{proof}
The results follow from the chain of equivalences:
\begin{equation*}
    \begin{split}
        LinCat_{\Hk}^{L,St} = & Fun^{Alg(\Sp)}(\un{\Hk},Pr^{L,St})        
      \\ \simeq  & \mathrm{LMod}_{\un{\Hk}}(Fun(*,Pr^{L,St})) \\ \simeq & \mathrm{LMod}_{\Hk}(Pr^{L,St}).
    \end{split}
\end{equation*}
The first equivalence holds by the definition of $\infty$-category of $Alg(\Sp)$-enriched $\infty$-functors plus the fact that the space of objects of $\un{\Hk}$ is trivial, i.e., $\un{\Hk}_{ob}\simeq *$.
The latter follows from the well-known equivalence $ev_{*}:Fun(*,Pr^{L,St})\simeq Pr^{L,St}$. 
\end{proof}

In the rest of this subsection, we will find other equivalent definitions using different higher structures of $Pr^{L,St}$. For example, using the fact that it is left tensored over itself (\ref{eqTensorPrlsT}), we obtain the following result: 

\begin{proposition}
\label{propNewDefinition}
There exists an equivalence of $\infty$-categories 
\[ LinCat_{\Hk}^{L,St}=\mathrm{LMod}_{\Hk}(Pr^{L,St})\simeq \mathrm{LMod}_{\LModk}(Pr^{L,St}).  \]
\end{proposition}

\begin{proof}
We prove the equivalence by showing that both the $\infty$-categories are apexes of the same pullback square.

The structure of $Pr^{L,St}$ as a left $Alg(\Sp)$-module object of $Pr^{L}$ is the restriction of its structure as a left $Pr^{L,St}$-module via the monoidal $\infty$-functor $\mathrm{LMod}_{(-)}:Alg(\Sp)\to Pr^{L,St}$. This means that it is obtained by the following pullback square:

\[\begin{tikzcd}[ampersand replacement=\&]
	{\mathrm{LMod}(Pr^{L,St})} \& {\mathrm{LMod}(Pr^{L,St})} \\
	{Alg(\Sp)} \& {Pr^{L,St}.}
	\arrow[from=1-1, to=1-2]
	\arrow[from=1-2, to=2-2]
	\arrow[from=2-1, to=2-2]
	\arrow[from=1-1, to=2-1]
	\arrow["\lrcorner"{anchor=center, pos=0.125}, draw=none, from=1-1, to=2-2]
\end{tikzcd}\]

By definition, \cite[Example 3.15]{DoniHigherCategorical}, the $\infty$-category of left $\Hk$-module object of $Pr^{L,St}$, there is the following pullback square:

\[\begin{tikzcd}[ampersand replacement=\&]
	{\mathrm{LMod}_{\Hk}(Pr^{L,St})} \& {\mathrm{LMod}(Pr^{L,St})} \\
	{\mathbb{1}} \& {Alg(\Sp).}
	\arrow[from=1-2, to=2-2]
	\arrow["{\{\Hk\}}", from=2-1, to=2-2]
	\arrow[from=1-1, to=1-2]
	\arrow[from=1-1, to=2-1]
	\arrow["\lrcorner"{anchor=center, pos=0.125}, draw=none, from=1-1, to=2-2]
\end{tikzcd}\]

Now, pasting together the above pullback squares, we obtain the diagram:

\begin{equation}
\label{eqCompositionSquares}
\begin{tikzcd}[ampersand replacement=\&]
	{\mathrm{LMod}_{\Hk}(Pr^{L,St})} \& {\mathrm{LMod}(Pr^{L,St})} \& {\mathrm{LMod}(Pr^{L,St})} \\
	{\mathbb{1}} \& {Alg(\Sp)} \& {Pr^{L,St}.}
	\arrow[from=1-2, to=1-3]
	\arrow[from=1-3, to=2-3]
	\arrow["{\mathrm{LMod}_{(-)}}"', from=2-2, to=2-3]
	\arrow[from=1-2, to=2-2]
	\arrow["{\{\Hk\}}"', from=2-1, to=2-2]
	\arrow[from=1-1, to=2-1]
	\arrow[from=1-1, to=1-2]
	\arrow["\lrcorner"{anchor=center, pos=0.125}, draw=none, from=1-1, to=2-2]
	\arrow["\lrcorner"{anchor=center, pos=0.125}, draw=none, from=1-2, to=2-3]
\end{tikzcd}
\end{equation}

Using the composition and cancellation of pullback squares rule, \cite[Proposition 4.3.11]{RiehlElements}, we deduce that the composite rectangle is a pullback.

By definition, \cite[Example 3.15]{DoniHigherCategorical}, the $\infty$-category of left $\LModk$-module object of $Pr^{L,St}$, there is the following pullback square:

\begin{equation}
\label{eqSquareLModk}
\begin{tikzcd}[ampersand replacement=\&]
	{\mathrm{LMod}_{\LModk}(Pr^{L,St})} \&\& {\mathrm{LMod}(Pr^{L,St})} \\
	{\mathbb{1}} \&\& {Pr^{L,St}.}
	\arrow[from=1-1, to=2-1]
	\arrow[from=1-1, to=1-3]
	\arrow["{\{\LModk\}}"', from=2-1, to=2-3]
	\arrow[from=1-3, to=2-3]
	\arrow["\lrcorner"{anchor=center, pos=0.125, rotate=45}, draw=none, from=1-1, to=2-3]
\end{tikzcd}
\end{equation}

Since the external rectangle pullback square in (\ref{eqCompositionSquares}) and the pullback square in (\ref{eqSquareLModk}) are computed on the same diagram, their apexes are equivalent. Thus, 
\[\mathrm{LMod}_{\LModk}(Pr^{L,St})\simeq \mathrm{LMod}_{\Hk}(Pr^{L,St}).\]

\end{proof}

In the last proposition, we have almost recovered the definition of $k$-linear stable presentable $\infty$-category in \cite[Variant D.1.5.1.]{SAG}, but in our definition there is an extra subscription \enquote{st}. In the next proposition, we prove that the two definitions are equivalent.

\begin{proposition}
\label{propTogliereSt}
There is an equivalence:
\[ \mathrm{LMod}_{\LModk}(Pr^{L,St})\simeq \mathrm{LMod}_{\LModk}(Pr^{L}).  \]
\end{proposition}
\begin{proof}
We shall prove the proposition by demonstrating that the following chain of equivalences holds:
\begin{equation}
\label{eqLautaro}
    \begin{split}    
    \mathrm{LMod}_{\LModk}(Pr^{St}) \simeq & \mathrm{LMod}_{\LModk}(\mathrm{LMod}_{\Sp}(Pr^{L}))\simeq  \mathrm{LMod}_{\LModk}(\mathrm{RMod}_{\Sp}(Pr^{L})) \\ \simeq & \; _{\LModk} \mathrm{BMod}(Pr^{L})_{\Sp} 
    \simeq  \; \mathrm{LMod}_{\LModk\otimes^{L} \Sp }(Pr^{L}) \\ \simeq & \mathrm{LMod}_{\Sp(\LModk) }(Pr^{L})  
    \simeq  \mathrm{LMod}_{\LModk}(Pr^{L}). 
    \end{split}
\end{equation}

The first equivalence follows from the monoidal equivalence \[Pr^{L,St}\simeq \mathrm{Mod}_{Sp}(Pr^{L}),\] see (\ref{eqPrstModSpPrL}).

In (\ref{eqStabilizatorTensor}), we have shown that the stabilization $\infty$-functor is equivalent to tensorizing with $\Sp$, i.e., there is a natural equivalence $\Sp(-)\simeq -\otimes^{L}(\Sp):Pr^{L}\to Pr^{L,St}$. So the second-to-last one holds. 
Since $\LModk$ is presentable and stable, the last one also holds.
The second holds because $\Sp$ is a commutative algebra.
The third holds for \cite[Theorem 4.3.2.7.]{HA}.
The next one holds by the folding technique, see \cite[Corollary in \S 3.6 ]{HinichYon}. Indeed, it implies that the $\infty$-category of $\LModk$-$\Sp$-bimodule objects of $Pr^{L}$ and the $\infty$-category of left $\LModk\otimes^{L}\Sp^{rev}$-modules of $Pr^{L}$.
\end{proof}

Hence, we find Lurie's definition:

\begin{definition}[\cite{SAG} Variant D.1.5.1]
\label{defKStableLinearPresentablecategoryLurie}
We denote by $LinCat_k^{L,St}$ the $\infty$-category of \textit{$k$-linear stable $\infty$-categories}, the $\infty$-category $\mathrm{LMod}_{\LModk}(Pr^{L})$.  We will refer to the elements of $LinCat_k^{L,St}$
 as \textit{stable presentable $k$-linear
$\infty$-categories}.
\end{definition}



\subsection{Morita theory for $k$-linear stable $\infty$-categories }
\label{subsecMoriaTheoryklinear}

In \cite{blumberg2013universal}, the authors define a Morita theory for stable $\infty$-categories without considering a $k$-linearization, which in this paper we refer to as non-geometrical Morita theory. 

In this subsection, we will generalize their theory by incorporating a general $k$-linearization.

In \Cref{propequivalenceCatperfprlomegaalgspmod}, we proved that the $\infty$-categories $\Cat^{perf}$ and $Pr^{L,St}_{\omega}$ are equivalent as left $Alg(\Sp)$-module objects of $Pr^{L}$. In particular, this implies that their $k$-linearizations are equivalent.  

The idea presented in \Cref{eqdefinitionideaDefKLinear} leads us to the following two natural definitions.

\begin{definition}
\label{defKlinearomegastable}
We denote by $LinCat_{\Hk}^{St,\omega}$ the $\infty$-category $Fun^{Alg(\Sp)}(\un{\Hk}, Pr^{L,St}_{\omega})$ and call it the $\infty$-category \textit{of $k$-linear h-projective} $\infty$-categories. Moreover, we call \textit{$k$-linear h-projective $\infty$-categories} the objects of $LinCat_{\Hk}^{St,\omega}$. 
\end{definition}

\begin{definition}
\label{defKlinearperf}
We denote by $LinCat_{\Hk}^{perf}$ the $\infty$-category $Fun^{Alg(\Sp)}(\un{\Hk}, \Cat^{perf})$ and call it the $\infty$-category \textit{of $k$-linear perfect $\infty$-categories}. Moreover, we call \textit{$k$-linear perfect $\infty$-categories} the objects of $LinCat_{\Hk}^{perf}$. 
\end{definition}

Using the definition of $\infty$-category of enriched $\infty$-functors \cite[\S 5.4]{DoniHigherCategorical}, we can immediately give an alternative definition. This is analog to \Cref{propActionklinearcategory}.

\begin{proposition}
\label{propDescription}
There are two equivalences of $\infty$-categories:
\begin{align*}
    LinCat_{\Hk}^{St,\omega} & = Fun^{Alg(\Sp)}(\un{\Hk}, Pr^{L,St}_{\omega})\simeq \mathrm{LMod}_{\Hk}(Pr^{L,St}_{\omega}); \\
    LinCat_{\Hk}^{perf} & = Fun^{Alg(\Sp)}(\un{\Hk}, \Cat^{perf})\simeq \mathrm{LMod}_{\Hk}(\Cat^{perf}).
\end{align*}
\end{proposition}

\begin{proof}
The proof is similar to that of \Cref{propActionklinearcategory}.
\end{proof}

Using the left module structure in \Cref{rmkSptensorizationonCatPerfPrSt}, we also find the analogous definition of \Cref{propActionklinearcategory}.

\begin{proposition}
\label{propNewDefinitionProj}
There are four equivalences of $\infty$-categories: 
\begin{itemize}
\item[(1)]
\[ LinCat_{\Hk}^{St,\omega}\simeq \mathrm{LMod}_{\LModk^{\omega}}(Pr^{L,St}_{\omega});  \]
\item[(2)]
\[ LinCat_{\Hk}^{St,\omega}\simeq \mathrm{LMod}_{\LModk}(Pr^{L,St}_{\omega}).  \]
\item[(3)]
\[ LinCat_{\Hk}^{perf}\simeq \mathrm{LMod}_{\LModk^{\omega}}(\Cat^{perf}).  \]
\item[(4)]
\[ LinCat_{\Hk}^{perf}\simeq \mathrm{LMod}_{\LModk}(\Cat^{perf}).  \]
\end{itemize}
\end{proposition}

\begin{proof}
The proof is similar to that of \Cref{propNewDefinition}, except that we consider different left module structures.
Indeed, instead of considering the left $Pr^{L,St}$-module structure of $Pr^{L}$: 
\begin{itemize}
    \item  in $(1)$ we consider $Pr_{\omega}^{L,St}$ as a left $\Cat^{perf}$-module object of $Pr^{L}$;

    \item in $(2)$ we consider $Pr_{\omega}^{L,St}$ as a left $Pr^{L,St}_{\omega}$-module object of $Pr^{L}$;

    \item in $(3)$ we consider $\Cat^{perf}$ as a left $\Cat^{perf}$-module object of $Pr^{L}$;

    \item in $(4)$ we consider $\Cat^{perf}$ as a left $Pr^{L,St}_{\omega}$-module object of $Pr^{L}$.
\end{itemize}

\end{proof}

As mentioned before, the $\infty$-category of perfect $k$-linear $\infty$-categories and the $\infty$-category of h-project $k$-linear $\infty$-categories are equivalent. 

\begin{proposition}
\label{propkCatPerfequivalentPrLomega}
The $\infty$-category of perfect $k$-linear $\infty$-categories and the $\infty$-category of h-projective $k$-linear $\infty$-categories are equivalent:

\[  Ind(-):LinCat^{St,w}_{\Hk}\simeq LinCat^{perf}_{\Hk}: (-)^{\omega}.\]      
\end{proposition}

\begin{proof}
In \Cref{propequivalenceCatperfprlomegaalgspmod}, we proved that $\Cat^{perf}$ and $Pr^{L,St}_{\omega}$ are equivalent as left $Alg(\Sp)$-module objects of $Pr^L$, so their $\infty$-categories of left $\Hk$-module objects are equivalent:
\[\mathrm{LMod}_{\Hk}(Pr^{L,St}_{\omega})\simeq \mathrm{LMod}_{\Hk}(\Cat^{perf}).\]
\end{proof}

In particular, as a consequence of \Cref{propNewDefinitionProj} and \Cref{propkCatPerfequivalentPrLomega}, we find the result of Cohn \cite[Proposition 5.6]{Cohn}:

\begin{corollary}
There is an equivalence of $\infty$-categories:
\[\mathrm{LMod}_{\LModk^{\omega}}(\Cat^{perf})\simeq \mathrm{LMod}_{\LModk}(\Cat^{perf}).\]
\end{corollary}

Using the same idea as in \Cref{propTogliereSt}, we can eliminate the subscript $St$.

\begin{proposition}
\label{propeliminareStOmega}
   There exists an equivalence of $\infty$-categories:
    \[LinCat^{St,\omega}_{\Hk}\simeq \mathrm{LMod}_{\LModk}(Pr^{L}_{\omega}).\]
\end{proposition}

\begin{proof}
 In \Cref{propNewDefinitionProj} we proved that there is an equivalence \[LinCat_{\Hk}^{St,\omega}\simeq \mathrm{LMod}_{\LModk}(Pr^{L,St}_{\omega}),\] 
 so the proof is similar to that of \Cref{propTogliereSt}.
\end{proof}


In the final part of this subsection, we will explain why all the definitions of $LinCat^{perf}_{\Hk}$ and $LinCat^{St,\omega}$ can be rewritten using the $\infty$-categories $\Cat^{\Sp, perf}$ and $ \Cat^{\Sp,\omega}$ instead of the $\infty$-categories $\Cat^{perf}$ and $ Pr^L_{\omega}$.

To do this, we will demonstrate that $\Cat^{\Sp,perf}$ and $\Cat^{\Sp,\omega}$ have equivalent structures as left $Alg(\Sp)$-module objects of $Pr^{L}$; this structure will not be the restriction of the higher structures of the $\infty$-category of $\Sp$-enriched $\infty$-categories, \cite[\S 5.3]{DoniHigherCategorical}.

In the proof of \Cref{propCatSpPerfPrOmega}, we established that there are two adjoint equivalences: 
\begin{equation}
\label{eqMamma}
(-)_o:\Cat^{\Sp,\omega}\simeq  Pr^{L,St}_{\omega}:\chi ;
\end{equation}
and 
\begin{equation}
\label{eqMammaBella}
(-)_o:Cat^{\Sp,perf}\simeq \Cat^{perf}:(Y\chi\circ Ind(-)).
\end{equation}

In particular, $Cat^{\Sp,perf}$ and $\Cat^{\Sp,\omega}$ inherit a monoidal structure from $\Cat^{perf}$ and $ Pr^{L,St}_{\omega}$ and a structure as left $Alg(\Sp)$-module objects in $Pr^{L}$. Thus, the equivalence \eqref{eqimportant} becomes a monoidal equivalence and also an equivalence of left $Alg(\Sp)$-module objects in $Pr^{L}$.

\begin{remark}
\label{rmkMonoidalStructure}
Now we summarize the higher structures of $\Cat^{\Sp,perf}$ and $\Cat^{\Sp,\omega}$:
\begin{itemize}

    \item $\Cat^{\Sp,perf}$ is a presentably monoidal $\infty$-category:
    \[-\otimes^{\Sp,perf}-:\Cat^{\Sp,perf} \times \Cat^{\Sp,perf}\to \Cat^{\Sp,perf}:\]
    with the unit being the $\Sp$-enriched $\infty$-category of compact spectra denoted by $Y\chi\circ Ind(\Sp^{\omega})$;

    \item $\Cat^{\Sp,perf}$ is a left $Alg(\Sp)$-module object of $Pr^{L}$:
    \[-\;^{\Sp}\otimes^{\Sp,perf} -:Alg(\Sp) \times \Cat^{\Sp,perf}\to \Cat^{\Sp,perf}:\]

    \item $\Cat^{\Sp,\omega}$ is a presentably monoidal $\infty$-category:
    \[-\otimes^{\Sp,\omega}-:\Cat^{\Sp, \omega} \times \Cat^{\Sp, h-\omega}\to \Cat^{\Sp,\omega}:\]
     with the unit being the $\Sp$-enriched $\infty$-category of spectra denoted by $\chi\circ \Sp$;

    \item  $\Cat^{\Sp,\omega}$ is a left $Alg(\Sp)$-module object of $Pr^{L}$:
    \[-\;^{\Sp}\otimes^{\Sp,\omega}-:Alg(\Sp) \times \Cat^{\Sp, h\text{-}proj}\to \Cat^{\Sp,\omega}.\]

\end{itemize}    


\end{remark}

\begin{corollary}
\label{corCatSpperfcatSpomegaequivalencemonoidalandModule}
The equivalences \eqref{eqimportant}, \eqref{eqMamma}, and \eqref{eqMammaBella} extend to a monoidal equivalence and to an equivalence of left $Alg(\Sp)$-module objects in $Pr^{L}$, where $\Cat^{\Sp,perf}$ and $\Cat^{\Sp,\omega}$ have the structures in \Cref{rmkMonoidalStructure} and $\Cat^{perf}$ and $Pr^{L,St}_{\omega}$ have the structure in \Cref{rmkSptensorizationonCatPerfPrSt}.
  
\end{corollary}

As an immediate consequence of the above corollary, we obtain that all the definitions of $LinCat^{perf}_{\Hk}$ and $LinCat^{St,\omega}$ found in this subsection (\Cref{propNewDefinitionProj},  \Cref{propDescription}, \Cref{defKlinearomegastable}, and \Cref{defKlinearperf}) can be reformulated with the $\infty$-categories $\Cat^{\Sp, perf}$ and $ \Cat^{\Sp,\omega}$ instead of the $\infty$-categories $\Cat^{perf}$ and $ Pr^L_{\omega}$. For example, $LinCat^{perf}_{\Hk}\simeq \mathrm{LMod}_{\Hk}(\Cat^{\Sp,perf})\simeq \mathrm{LMod}_{\LModk}(\Cat^{\Sp,\omega})$.

\section{Results}
\label{secRacks}

In this section, we will compare the two approaches to linearization seen in \Cref{secDg} and \Cref{secklinearstableCat}.

In the first part, we will define a localization in the non-geometrical setting, which is the $\infty$-categorical analog of the localization $\Psi_{\text{perf}}$ in \cite[Theorem 4.23]{blumberg2013universal}. 
From \Cref{remarkCAZZOCAZZOCAZZO}, we know that there is a localization:
\begin{equation}
\label{eqLOcaCAZZOCAZZOCAZZO}
\begin{tikzcd}[ampersand replacement=\&]
	{\Cat^{\Sp}} \&\& {\Cat^{\Sp,St}\simeq\Cat^{Ex}}
	\arrow[""{name=0, anchor=center, inner sep=0},"{\mathscr{S}(-)_{o}}", shift left=2, from=1-1, to=1-3]
	\arrow[""{name=1, anchor=center, inner sep=0}, shift left=2, hook', from=1-3, to=1-1]
	\arrow["\dashv"{anchor=center, rotate=-90}, draw=none, from=0, to=1]
\end{tikzcd}
\end{equation}

In \cite[Lemma 2.20]{blumberg2013universal}, the authors prove that there is a localization: 
\begin{equation}
 \label{eqFIGAFIGAFIGA}
\begin{tikzcd}[ampersand replacement=\&]
	{\Cat^{Ex}} \&\& {\Cat^{perf}}
	\arrow[""{name=0, anchor=center, inner sep=0}, "{Idem(-)}", shift left=2, from=1-1, to=1-3]
	\arrow[""{name=1, anchor=center, inner sep=0}, shift left=2, hook', from=1-3, to=1-1]
	\arrow["\dashv"{anchor=center, rotate=-90}, draw=none, from=0, to=1]
\end{tikzcd}
\end{equation}
whose left adjoint is the $\infty$-functor in \Cref{exIdempo}, which can also be described on the objects as $Idem(\C) \simeq (Ind(\C))^{\omega}$ as we mentioned in \Cref{propIdemComplCompact}.

Moreover, we know that there is an equivalence $Ind(-):\Cat^{\text{perf}}\simeq Pr^{L,St}_{\omega}$, see \eqref{eqCatPerfPrStLOmega}.

Putting together the above we obtain the Morita localization in the non-geometrical case. 

\begin{proposition}
\label{proprimaLocalizazioneZCaso} 
There are two localizations of $\infty$-categories whose right adjoint is lax monoidal:   

\[\begin{tikzcd}[ampersand replacement=\&]
	{\Cat^{\Sp}} \&\& {{Pr^{L,St}_{\omega}};} \\
	{\Cat^{\Sp}} \&\& {\Cat^{perf}.}
	\arrow[""{name=0, anchor=center, inner sep=0}, "{Ind(Idem(\mathscr{S}(-)_{o}))}", shift left=2, from=1-1, to=1-3]
	\arrow[""{name=1, anchor=center, inner sep=0}, shift left=2, hook', from=1-3, to=1-1]
	\arrow[""{name=2, anchor=center, inner sep=0}, "{Idem(\mathscr{S}(-)_{o})}", shift left=2, from=2-1, to=2-3]
	\arrow[""{name=3, anchor=center, inner sep=0}, shift left=2, hook', from=2-3, to=2-1]
	\arrow["\dashv"{anchor=center, rotate=-90}, draw=none, from=0, to=1]
	\arrow["\dashv"{anchor=center, rotate=-90}, draw=none, from=2, to=3]
\end{tikzcd}\]

\end{proposition}

\begin{proof}
The second localization has been proved before the proposition was stated.
The first is the second postcomposed with the equivalence (\ref{eqCatPerfPrStLOmega}). 
\end{proof}

Using \Cref{corCatSpperfcatSpomegaequivalencemonoidalandModule}, we obtain the following corollary.

\begin{corollary}
\label{corFanculo}
There exists two localization of $\infty$-categories:
\[\begin{tikzcd}[ampersand replacement=\&]
	{\Cat^{\Sp}} \&\& {\Cat^{\Sp, \omega};} \\
	{\Cat^{\Sp}} \&\& {\Cat^{\Sp, perf}.}
	\arrow[""{name=0, anchor=center, inner sep=0}, "{\mathfrak{T}}", shift left=2, from=1-1, to=1-3]
	\arrow[""{name=1, anchor=center, inner sep=0}, "\iota", shift left=2, hook', from=1-3, to=1-1]
	\arrow[""{name=2, anchor=center, inner sep=0}, "{\mathfrak{M}}", shift left=2, from=2-1, to=2-3]
	\arrow[""{name=3, anchor=center, inner sep=0}, "\iota", shift left=2, hook', from=2-3, to=2-1]
	\arrow["\dashv"{anchor=center, rotate=-90}, draw=none, from=0, to=1]
	\arrow["\dashv"{anchor=center, rotate=-90}, draw=none, from=2, to=3]
\end{tikzcd}\]
\end{corollary}

Now we want to understand the essential image of the right adjoint in \Cref{proprimaLocalizazioneZCaso} and \Cref{corFanculo}, and we will describe it more explicitly.

In \cite[Theorem 1.2]{heine2023equivalence}, the author proves that there exists a fully faithful $\infty$-functor 
\[\chi:P^{L,St}\simeq \mathrm{LMod}_{\Sp}(Pr^{L})\to \widehat{\Cat}^{\Sp}\]
 whose essential image is the $\infty$-category $\Cat^{\Sp,L,St}$, which we defined in \Cref{defcategorypresentablespectral}. The first equivalence is \eqref{eqPrstModSpPrL}.

First, we restrict  $\chi$ along the inclusion $P^{L,St}_{\omega}\subseteq Pr^{L,St}_{\omega}$ and we obtain a fully faithful $\infty$-functor:
\begin{equation}
\label{equationCAZZAROLissima}
\chi_{|Pr^{L,St}_{\omega}}:Pr^{L,St}_{\omega}\to \widehat{Cat_{\infty}}^{\Sp}.
\end{equation}
The essential image of $\chi_{|Pr^{L,St}_{\omega}}$ is the $\infty$-category of h-projective spectral $\infty$-categories $\Cat^{\Sp,\omega}$, which we defined in \Cref{defSpectralCategoryofhprojectomega}. This is not the correct $\infty$-functor in our localization because the objects of $\Cat^{\Sp,\omega}$ are not necessarily small.

The Morita theory for spectral $\infty$-categories solves this problem.
In \Cref{propCatSpPerfPrOmega}, we proved that $\Cat^{\Sp,\omega}$ is equivalent to the $\infty$-category of perfect spectral $\infty$-categories $\infty$-categories $\Cat^{\Sp,perf}$, so we can change the $\infty$-functor $\chi_{|Pr^{L,St}_{\omega}}$ so that its essential image is $\Cat^{\Sp,perf}$. Now we have overcome the problem because the objects of $\Cat^{\Sp,perf}$ are small by definition. 
So the $\infty$-functor $\chi_{|Pr^{L,St}_{\omega}}$ factors through the inclusion $\Cat^{\Sp}\subseteq \widehat{\Cat}^{\Sp}$, and we obtain the correct right adjoint in our localization, the lax monoidal fully faithful: 
\begin{equation}
\label{equationCAZZAROLA}
\chi_{|Pr^{L,St}_{\omega}}:Pr^{L,St}_{\omega}\to \Cat^{\Sp},
\end{equation}
whose essential image is $\Cat^{\Sp,perf}$. 
Note that we use the same notation for the $\infty$-functors (\ref{equationCAZZAROLA}) and (\ref{equationCAZZAROLissima}) because they are essentially the same $\infty$-functor.
The $\infty$-functor (\ref{equationCAZZAROLA}) is the right adjoint in \Cref{proprimaLocalizazioneZCaso} and \Cref{corFanculo}. 
 
Above we have proved the following corollary.

\begin{corollary}
There exists an equivalence of $\infty$-categories:
\begin{equation}
\label{eqMainTh1}
\Cat^{\Sp,perf}\simeq Pr^{L,St}_{\omega}\simeq \Cat^{perf}.
\end{equation}
\end{corollary}

\begin{proof}
The $\infty$-functor (\ref{equationCAZZAROLA}) gives the first equivalence.
While the second is (\ref{eqCatPerfPrStLOmega}).
\end{proof}

Thanks to \Cref{CorCatLModPerfCatLModhprojequivalence} and \Cref{rmksuperimportant}, the same construction above holds if we replace the presentably monoidal $\infty$-category $\Sp$ by the presentably monoidal $\infty$-category $\LModk$.


For the sake of completeness, we will sketch the argument once more.

In \cite[Theorem 1.2]{heine2023equivalence}, the author proves that there exists a fully faithful $\infty$-functor 
\[\chi: \mathrm{LMod}_{\LModk}(Pr^{L,St})\simeq \mathrm{LMod}_{\LModk}(Pr^{L})\to \widehat{\Cat}^{\LModk},\]
whose essential image is the $\infty$-category $\Cat^{\LModk,L,St}$, which we define in \Cref{defdgCategoryofhproject}. Instead, the first equivalence is (\ref{propTogliereSt}).

We restrict $\chi$ along the inclusion $\mathrm{LMod}_{\LModk}(P^{L,St}_{\omega})\subseteq  Pr^{L,St}_{\omega}$ and we obtain a fully faithful $\infty$-functor:
\begin{equation}
\label{equationCAZZAROLissimaLModk}
\chi_{|\mathrm{LMod}_{\LModk}(P^{L,St}_{\omega})}:\mathrm{LMod}_{\LModk}(Pr^{L,St}_{\omega})\to \widehat{Cat_{\infty}}^{\LModk}.
\end{equation}

The essential image of $\chi_{|\mathrm{LMod}_{\LModk}(P^{L,St}_{\omega})}$ is the $\infty$-category $\Cat^{\LModk,\omega}$, which we defined in \Cref{defdgCategoryofhomega}. To use it in our last proposition we need the $\infty$-functor (\ref{equationCAZZAROLissimaLModk}) to target the $\infty$-category of small $\LModk$-enriched $\infty$-categories.
However, the objects $\Cat^{\LModk,\omega}$ are not necessarily small, so this $\infty$-functor is not the right one.
The Morita theory for dg-categories resolves this problem. 
As shown in \Cref{CorCatLModPerfCatLModhprojequivalence}, the $\infty$-category $\Cat^{\LModk,\omega}$ is equivalent to the $\infty$-category of perfect spectral $\infty$-categories $\Cat^{\LModk,perf}$, allowing us to adjust the $\infty$-functor $\chi_{|\mathrm{LMod}_{\LModk}(P^{L,St}_{\omega})}$ so that its essential image is $\Cat^{\LModk,perf}$. This resolves the issue as the objects of $\Cat^{\LModk,perf}$ are small by definition. 

We find that the $\infty$-functor $\chi_{|\mathrm{LMod}_{\LModk}(P^{L,St}_{\omega})}$ factors through the inclusion $\Cat^{\LModk}\subseteq \widehat{\Cat}^{\LModk}$, allowing us to obtain the fully faithful $\infty$-functor: 
\begin{equation}
\label{equationCAZZAROLALMod}
\chi_{|\mathrm{LMod}_{\LModk}(P^{L,St}_{\omega})}:\mathrm{LMod}_{\LModk}(P^{L,St}_{\omega})\to Cat^{\LModk},
\end{equation}
whose essential image is $\Cat^{\LModk,perf}$.
Note that we use the same notation for the $\infty$-functors (\ref{equationCAZZAROLALMod}) and (\ref{equationCAZZAROLissimaLModk}) because they are essentially the same $\infty$-functor.

Above, we have proven the following result.

\begin{proposition}
\label{Useful}
Let $k$ be a discrete $\E$-ring. There exists an equivalence of $\infty$-categories:
\begin{equation}
\label{eqMainTh3}
\Cat^{\LModk,perf}\simeq \mathrm{LMod}_{\LModk}(P^{L,St}_{\omega}).
\end{equation}
\end{proposition}

Finally, we can state the Main Theorem of this paper: \enquote{idempotent complete pretriangulated dg-categories over $k$ are stable $k$-linear idempotent complete $\infty$-categories}.

\begin{theorem}[Main Theorem]
\label{thMaintheorem}
Let $k$ be a discrete $\E$-ring. There exists an equivalence of $\infty$-categories:
\begin{equation}
\label{eqMainTh2}
\Cat^{\D(k),perf}\simeq LinCat^{perf}_{\Hk}.
\end{equation}
\end{theorem}
\begin{proof}
This follows easily from \Cref{Useful}, \Cref{propEqDkenrLModkenric}, and \Cref{propNewDefinitionProj}.
\end{proof}

Now, we define what it means for two dg-categories to be Morita equivalent.
We denote by $\mathscr{U}$ the following composition of  $\infty$-functors:

\begin{equation}
\label{eqFunctorThatDefine}
\mathscr{U}:\Cat^{\D(k)}\simeq \Cat^{\LModk}\simeq \mathrm{LMod}_{\Hk}(\Cat^{\Sp})\xrightarrow{U}\Cat^{\Sp}\xrightarrow{\mathfrak{M}} \Cat^{\Sp,perf},
\end{equation}
we found the first equivalence in \Cref{propEqDkenrLModkenric}, the second equivalence is \cite[Corollary 6.9]{DoniHigherCategorical}, $U$ is the canonically forgetful $\infty$-functor, and $\mathfrak{M}$ is the $\infty$-functor in \Cref{corFanculo}.

\begin{definition}
\label{defMoritaEquivalence}
Let $\A$ and $\B$ be two objects of $\Cat^{\D(k)}$, and $F:\A\to \B$ be an arrow in $\Cat^{\D(k)}$ between them.
We say that $F$ is a \textit{Morita equivalence} if the arrow $\mathscr{U}(F):\mathscr{U}(\A)\to\mathscr{U}(\B) $ is an equivalence in $\Cat^{\Sp,perf}$.
\end{definition}

Let $\mathscr{M}_k$ be the collection of Morita equivalences of $\Cat^{\D(k)}$. Since $\Cat^{\D(k)}$ is a presentable $\infty$-category, we can localize with respect to the strongly saturated
class of morphisms generated by $\mathscr{M}_k$ \cite[Proposition 5.5.4.15]{HTT}, and we obtain a localization

\begin{equation}
\label{eqLocalizazioneFinalmente}
\begin{tikzcd}[ampersand replacement=\&]
	{\Cat^{\D(k)}} \&\& {\Cat^{\D(k)}[\mathscr{M}_k^{-1}]}
	\arrow[""{name=0, anchor=center, inner sep=0}, "{\mathfrak{i}}", shift left=3, hook', from=1-3, to=1-1]
	\arrow[""{name=1, anchor=center, inner sep=0}, "{\mathfrak{L}}", shift left=3, from=1-1, to=1-3]
	\arrow["\dashv"{anchor=center, rotate=-90}, draw=none, from=1, to=0]
\end{tikzcd}    
\end{equation}

In the last result of this paper, we proved that the $\infty$-category $\Cat^{\LModk}[\mathscr{M}_k^{-1}]$ is equivalent to the $\infty$-category $LinCat_{\Hk}^{perf}$.

\begin{proposition}
\label{propklinearLocalMorita}
There exists a localization
\begin{equation}
\label{basta}
\begin{tikzcd}[ampersand replacement=\&]
	{\Cat^{\D(k)}} \&\& {LinCat^{perf}_{\Hk},}
	\arrow[""{name=0, anchor=center, inner sep=0}, "{\mathfrak{i}}", shift left=3, hook', from=1-3, to=1-1]
	\arrow[""{name=1, anchor=center, inner sep=0}, "{\mathfrak{M}_k}", shift left=3, from=1-1, to=1-3]
	\arrow["\dashv"{anchor=center, rotate=-90}, draw=none, from=1, to=0]
\end{tikzcd}
\end{equation}
\end{proposition}

\begin{proof}
For \Cref{Useful}, \Cref{propEqDkenrLModkenric}, and \Cref{propNewDefinitionProj}, to prove the results, it is enough to prove that $\Cat^{\LModk}[\mathscr{M}_k^{-1}]$ is equivalent to $\Cat^{\LModk, perf}$. Since \eqref{eqLocalizazioneFinalmente} is a localization of a presentable $\infty$-category, the above sentence is proved if we show that a dg-category is perfect if and only if it is a $\mathscr{M}_k$-local object.
Since $\Cat^{\Sp}$ and $\Cat^{\Sp,perf}\simeq \Cat^{perf}$ are presentable $\infty$-categories the localization in \Cref{corFanculo} is characterized by its class of local morphisms, \cite[Proposition 5.5.4.2]{HTT}, we denote this class by $\mathscr{M}$.

We start with a claim:
\begin{displayquote}
(*) if $f:\B\to \C$ belongs to $\mathscr{M}$, then $\Hk\otimes f: \Hk\otimes^{\Sp} \B\to \Hk\otimes^{\Sp}\C$ belongs to $\mathscr{M}$.
\end{displayquote}
Let's assume for a moment that the claim holds.

Let
\begin{equation}
\label{eqAggiunzioneExtentionRestriction}    
F: \Cat^{\Sp}\rightleftarrows \mathrm{LMod}_{\Hk}(\Cat^{\Sp}):U
\end{equation}
be the  extension-restriction of scalar adjunction, see \cite[\S 4.2]{DoniHigherCategorical}.

By definition of $\mathscr{M}_{k}$, \eqref{eqFunctorThatDefine}, $U$ sends arrows in $\mathscr{M}_{k}$ into arrows in $\mathscr{M}$.

We show that $F$ sends arrows in $\mathscr{M}$ into arrows in $\mathscr{M}_k$.

Let $f:\B\to \C$ be an element in $\mathscr{M}$ and we need to show that the $\infty$-functor \eqref{eqFunctorThatDefine} sends $F(f): F(\B)\to F(\C)$ into an isomorphism.  
In short, we need to prove that $UF(f): UF(\B)\to UF(\C)$ belongs to $\mathscr{M}$. This is a consequence of claim $(*)$, because the $\infty$-functor $UF(-)$ is equivalent to the $\infty$-functor $\Hk\otimes^{\Sp}-:\Cat^{\Sp}\to \Cat^{\Sp}$, see \cite[\S 4.2]{DoniHigherCategorical}.

Let $\A\in \mathrm{LMod}_{\Hk}(\Cat^{\Sp})$ be an $\mathscr{M}_{k}$-local object and we want to prove that it also belongs to $\mathrm{LMod}(\Cat^{\Sp, perf})$; i.e. that $U(\A)$ is a perfect spectral $\infty$-category or equivalently that it is a $\mathscr{M}$-local object of $\Cat^{\Sp}$. 
So, let $f:\B\to \C$ be an arrows in $\mathscr{M}$, we need to prove that the composition morphism 
\[f^{*}:=(-)\circ f:\Hom_{\Cat^{\Sp}}(\C,U(\A))\to \Hom_{\Cat^{\Sp}}(\B,U(\A)),\]
is an equivalence of $\infty$-groupoids (or spaces); as usual, we denote the precomposition with an arrow $f$ by $f^{*}$.
Using the universal property of the adjunction \eqref{eqAggiunzioneExtentionRestriction}, it is equivalent to prove that the composition morphism   
\[F(f)^{*}:\Hom_{\mathrm{LMod}_{\Hk}(\Cat^{\Sp})}(F(\C),\A))\to \Hom_{\mathrm{LMod}_{\Hk}(\Cat^{\Sp})}(F(\B),\A),\]
is an equivalence of $\infty$-groupoids.
But this is holds because $\A$ is a $\mathscr{M}_{k}$-local object and $F$ sends arrows in $\mathscr{M}$ into arrows in $\mathscr{M}_k$.

Vice versa, let $\A\in \mathrm{LMod}_{\Hk}(\Cat^{\Sp})$ be a perfect dg-category we want to prove that $\A$ is a $\mathscr{M}_{k}$-local object of $\mathrm{LMod}_{\Hk}(\Cat^{\Sp})$. 

Let $f:\B\to \C$ be an arrow in $\mathscr{M}_{k}$ and we want to proves that the composition morphism
\begin{equation}
\label{eqBhoNonsoCheNotazioneDare}
f^{*}:\Hom_{\mathrm{LMod}_{\Hk}(\Cat^{\Sp})}(\C,\A)\to \Hom_{\mathrm{LMod}_{\Hk}(\Cat^{\Sp})}(\B,\A),
\end{equation}
is an equivalence of $\infty$-groupoids.

The morphism \eqref{eqBhoNonsoCheNotazioneDare} is the limit 
\begin{equation}
\label{eqINutileMaLAMetto}
    lim \Hom_{\Cat^{\Sp}}((\otimes^{\Sp,n}\Hk)\otimes^{\Sp}U(\B),U(\A))\to lim (\Hom_{\Cat^{\Sp}}((\otimes^{\Sp,n}\Hk)\otimes^{\Sp}U(\B),U(\A)))
\end{equation}

of the morphism between cosimplicial objects:
\[\begin{tikzcd}[ampersand replacement=\&]
	{\Hom_{\Cat^{\Sp}}(U(\C),U(\A))} \&\&\& {\Hom_{\Cat^{\Sp}}(U(\B),U(\A))} \\
	{\Hom_{\Cat^{\Sp}}(\Hk\otimes^{\Sp}U(\C),U(\A))} \&\&\& {\Hom_{\Cat^{\Sp}}(\Hk\otimes^{\Sp}U(\B),U(\A))} \\
	{\Hom_{\Cat^{\Sp}}(\Hk\otimes^{\Sp}\Hk\otimes^{\Sp}U(\C),U(\A))} \&\&\& {\Hom_{\Cat^{\Sp}}(\Hk\otimes^{\Sp}U(\B),U(\A))} \\
	{\Hom_{\Cat^{\Sp}}((\otimes^{\Sp,n}\Hk)\otimes^{\Sp}U(\C),U(\A))} \&\&\& {\Hom_{\Cat^{\Sp}}((\otimes^{\Sp,n}\Hk)\otimes^{\Sp}U(\B),U(\A))} \\
	{} \&\&\& {}
	\arrow[shift left=3, from=2-1, to=1-1]
	\arrow["{(\Hk\otimes^{\Sp}f)^{*}}", from=2-1, to=2-4]
	\arrow["{f^{*}}", from=1-1, to=1-4]
	\arrow[shift left=3, from=2-4, to=1-4]
	\arrow[shift right=3, from=2-4, to=1-4]
	\arrow[shift right=3, from=2-1, to=1-1]
	\arrow["{(\Hk\otimes^{\Sp}f)^{*}}", from=3-1, to=3-4]
	\arrow[shift left=3, from=3-4, to=2-4]
	\arrow[shift right=3, from=3-4, to=2-4]
	\arrow[from=3-4, to=2-4]
	\arrow[shift left=3, from=3-1, to=2-1]
	\arrow[shift right=3, from=3-1, to=2-1]
	\arrow[from=3-1, to=2-1]
	\arrow["{((\otimes^{\Sp,n}\Hk)\otimes^{\Sp}f)^{*}}", from=4-1, to=4-4]
	\arrow[shift left, dashed, no head, from=4-1, to=3-1]
	\arrow[dashed, no head, from=4-1, to=5-1]
	\arrow[dashed, no head, from=5-4, to=4-4]
	\arrow[shift right, dashed, no head, from=4-4, to=3-4]
\end{tikzcd}\]
here, we denoted by $\otimes^{\Sp,n}\Hk$ the object 
$\underbrace{\Hk\otimes^{\Sp}\dots\otimes^{\Sp}\Hk}_{n\text{-}times}$.
The horizontal morphisms in the above diagram between cosimplicial objects are equivalences of $\infty$-groupoids: $U(\A)$ is a $\mathscr{M}$-local object of $\Cat^{\Sp}$ and, for each $n\in\N$, the \[(\otimes^{\Sp,n}\Hk)\otimes^{\Sp}f:(\otimes^{\Sp,n}\Hk)\otimes^{\Sp}U(\B)\to (\otimes^{\Sp,n}\Hk)\otimes^{\Sp}U(\C) \]
belongs to $\mathscr{M}$ as a consequence of the claim $(*)$ and the fact that $U$ sends arrows in $\mathscr{M}_{k}$ into arrows in $\mathscr{M}$.
This implies that \eqref{eqINutileMaLAMetto} and \eqref{eqBhoNonsoCheNotazioneDare} are equivalences.

To conclude we are left only to prove claim $(*)$. 
Let $f:\B\to \C$ be an arrow in $\mathscr{\W}$ and let $\A$ be a $\mathscr{M}$-local object.
We want to prove that 
\[(\Hk\otimes^{\Sp} f)^{*}:\Hom_{\Cat^{\Sp}}(\Hk\otimes^{\Sp}\C,\A)\to \Hom_{\Cat^{\Sp}}(\Hk\otimes^{\Sp}\B,\A),\]
is an equivalence of $\infty$-groupoids.
Since $(\Cat^{\Sp},\otimes^{\Sp})$ has an internal-hom $\Cat^{\Sp}(-,-)$, it is the same prove that 
\begin{equation}
f^{*}:\Hom_{\Cat^{\Sp}}(\C,\Cat^{\Sp}(\un{\Hk},\A))\to \Hom_{\Cat^{\Sp}}(\B,\Cat^{\Sp}(\un{\Hk},\A)),
\end{equation}
is an equivalence of $\infty$-groupoids.
To archive this, it is sufficient to prove that $\Cat^{\Sp}(\un{\Hk},\A)$ is a $\mathscr{M}$-local object of $\Cat^{\Sp}$, i.e. it is perfect spectral $\infty$-category.
So we need to prove that its underlying $\infty$-category $\Cat^{\Sp}(\un{\Hk},\A)_{o}$ is perfect.

The following chain of equivalences holds:
\begin{equation}
\begin{split}
\Cat^{\Sp}(\un{\Hk},\A)_{o}\simeq & Fun^{\Sp}(\un{\Hk},\A)  \simeq Fun_{\Cat^{Ex}}(\mathscr{S}(\un{\Hk})_{o},\A_{o}) \\ \simeq & Fun_{\Cat^{perf}}(Idem(\mathscr{S}(\un{\Hk})_{o}),\A_o)\\ \simeq & Fun^{Ex}(Idem(\mathscr{S}(\un{\Hk})_{o}),\A_o);
\end{split}
\end{equation}
the first equivalence is by definition \cite[\S 5.4]{DoniHigherCategorical}, the second one is \cite[Corollary 9.20]{heine2023equivalence}, 
the third is by universal property of Idem, \Cref{exIdempo}, and the last one follows because $\Cat^{perf}$ is a full $\infty$-subcategory of $\Cat^{Ex}$.

In \cite[\S 4]{hoyois2017higher}, the authors prove that $Fun^{Ex}(-,-)$ is the internal-hom of $\Cat^{perf}$, so $\Cat^{\Sp}(\un{\Hk},\A)_{o}$ is perfect.
\end{proof}

Note that the functor $\mathfrak{i}$ in \eqref{basta} is equivalent to \eqref{equationCAZZAROLALMod}.

\vskip 1cm


\bibliographystyle{unsrt}
\bibliography{CategoricalArticle}

\end{document}